\documentclass[12pt,a4paper]{article}
\usepackage{fleqn,amsfonts,amssymb,amsmath,amscd,epsfig,version,theorem} %
\usepackage{fancyhdr,enumerate,epsfig}

\usepackage{yhmath}
\usepackage{hyperref}

\numberwithin{equation}{section}
\numberwithin{figure}{section}
\theoremstyle{change}
\newtheorem{theorem}{Theorem} [section]
\newtheorem {lemma}[theorem]{Lemma}

\newtheorem {corollary}[theorem]{Corollary}

\newtheorem {defi}[theorem]{Definition}
{\theorembodyfont{\normalfont}
{\theorembodyfont{\normalfont}\newtheorem {example}[theorem]{Example}
{\theorembodyfont{\normalfont}
{\theorembodyfont{\normalfont}\newtheorem {remark}[theorem]{Remark}
{\theorembodyfont{\normalfont}\newtheorem {remarks}[theorem]{Remarks}
\newcommand{\beq}{\begin{equation}}
\newcommand{\eeq}{\end{equation}}
\newcommand{\Leq}[1]{\label{#1}\end{equation}}
\newcommand{\beqn}{\begin{eqnarray}}
\newcommand{\eeqn}{\end{eqnarray}}
\newcommand{\beqno}{\begin{eqnarray*}}
\newcommand{\eeqno}{\end{eqnarray*}}
\newcommand{\es}{\emptyset}
\renewcommand {\l}{\left}
\newcommand {\ri}{\right}
\newcommand {\vv}{\varphi}     

\newcommand {\LA}{\left\langle}
\newcommand {\RA}{\right\rangle}

\newcommand {\eh}{{\textstyle \frac{1}{2}}}

\newcommand {\bN}{{\mathbb N}}
\newcommand {\bR}{{\mathbb R}}

\newcommand {\bS}{{\mathbb S}}

\newcommand{\idty}{{\rm 1\mskip-4mu l}} 
\newcommand{\cC}{{\cal C}} %
\newcommand{\cD}{{\cal D}} %

\newcommand{\cM}{{\cal M}}
\newcommand{\cP}{{\cal P}}

\newcommand{\bem}{\l(\! \begin{array}}
\newcommand{\eem}{\end{array}\!\ri)}
\newcommand{\bsm}{\left(\begin{smallmatrix}} 
\newcommand{\esm}{\end{smallmatrix}\right)}  
\newcommand{\NN}{\nonumber}

\newcommand{\hP}{{\widehat{P}\,}}

\newcommand{\qmbox}[1]{\quad\mbox{#1}\quad}
\renewcommand {\max}{{{\rm max}}}

\newcommand{\hSE}{{\widehat{\Sigma}_E}}
\newcommand{\wSE}{{\wideparen{\Sigma}_E}}

\def\XXint#1#2#3{{\setbox0=\hbox{$#1{#2#3}{\int}$}
\vcenter{\hbox{$#2#3$}}\kern-.5\wd0}}

\newcommand{\qe}{q^E} 
\newcommand{\pe}{p^E} 
\newcommand{\qi}{q^I} 
\newcommand{\ppi}{p^I} 
\newcommand{\De}{\Delta^E} 
\newcommand{\Di}{\Delta^I} 
\usepackage{color}

\begin{document}

\title{Compactification of the energy surfaces\\ for $n$ bodies}
\author{
Andreas Knauf\thanks{
Department of Mathematics,
Friedrich-Alexander-University Erlangen-N{\"u}rnberg,
Cauerstr.\ 11, D-91058 Erlangen, 
Germany, \texttt{knauf@math.fau.de}} 
\and
Richard Montgomery\thanks{Mathematics Department,
UC Santa Cruz, 
4111 McHenry, 
Santa Cruz, CA 95064, USA, \texttt{rmont@ucsc.edu}}}

\date{\today}
\maketitle
\begin{abstract}
For $n$ bodies moving in Euclidean $d$--space under the influence of a   
homogeneous pair interaction we 
compactify every center of mass energy surface, 
obtaining a  $(2d(n -1)-1)$--dimensional  manifold with corners in the sense of Melrose. 
After a time change, the flow on this manifold is globally defined 
and non-trivial on the boundary.
\end{abstract}
\tableofcontents
%
\section{Introduction and Results}
%
The $n$--body problem of celestial mechanics is described by the Hamiltonian\label{Hamiltonian}
\beq 
H(q,p) = K(p) - U(q) \mbox{ , with } 
U(q_1,\ldots,q_n) = \sum_{1\le i<j\le n} \frac{Z_{i,j}}{\|q_i-q_j\|^\alpha} \, ,
\Leq{Ham:Pot}
kinetic energy $K(p) = \sum_{i=1}^n \frac{\|p_i\|^2}{2m_i}$, $\alpha:=1$ and
$Z_{i,j}:= G m_i m_j$ with $G$ the gravitational constant. 
We call $U$ the `potential' although it is actually the negative of the potential.\\
We work in the center of mass system.  
Due to collisions and escape to spatial infinity its 
energy surfaces $\hSE = H^{-1}(E)$ are non-compact. Our goal here is to 
compactify the energy surfaces
by adding boundary pieces in such a way that the $n$--body flow extends 
nontrivially to the added boundaries. 
\begin{theorem} \quad\\ 
Assume $Z_{i,j} > 0$ and $0 < \alpha < 2$
in the `potential' $U$ of equation \eqref{Ham:Pot}.  Then the   constructions of sections 
\ref{sub:blow:conf} and \ref{sec:blow:energy} yield, for each energy $E$,  a compact manifold
with corners   $\wSE$ whose interior is diffeomorphic to the usual energy level set 
$\hSE$ away from the Hill boundary $\{U = - E\}$ 
and on which the usual   $n$--body flow, after a time reparameterization, extends
continuously to  a non-trivial $C^1$--flow on the added boundaries.
\end{theorem} 
For the definition of a manifold with corners and a brief introduction to these spaces 
see Appendix \ref{sect:appendix}.

We prove the theorem by combining two  well-established techniques 
from the $n$--body world with a  technique primarily used in linear PDEs.   

The first  technique is that of adding sub-manifolds to the standard  phase space 
in order to better understand limiting behaviours. For collisions
we add  collision manifolds by  blow up,  following the ideas of McGehee.  
See \cite{McG1,McG2,ElB,Dev,LS, LI} and references cited therein.
For  widely separated particles with interparticle distances diverging 
we add a  manifold at infinity.  See  \cite{DMMY}. 

The second technique,  originally developed for the quantum 
$n$--body problem, is that of {\sc Graf} partitions (see \cite{Gr}).
These partitions help us to label the strata of $\wSE$ and  understand 
its topology and combinatorics.  

The third technique is  Melrose's iterative 
real blow-up \cite{Me} developed for PDEs.  This iterative blow-up 
provides the appropriate language to  get to the final result. 
In the last section we deduce the  topology of the total blow up in 
Theorem~\ref{thm:boundary:hom}. 

The famous regularizations of binary collisions due to {\sc Levi-Civita} \cite{LC} or
{\sc Moser} \cite{Mo} (or  any other method)  
allow us to analytically pass through binary collisions 
but are only available when $\alpha = 2(1-1/k)$ for $k\in\bN$. 
(See \cite[Theorem 7.1]{McG2}).
Wanting a unified picture of compactification  valid  for all values of $\alpha$, 
or at least values
in the interval $(0,2)$, we have had to treat binary collisions in  the same way 
as we have treated triple or higher collisions: 
by slowing down the flow and attaching boundaries corresponding to these collisions.

We will  work in the center of mass configuration space 
\[ \textstyle
M = \big\{q=(q_1,\ldots,q_n)\in\bR^{nd} \mid \sum_{i=1}^n m_i q_i = 0\big\} .\] 
The (negative) potential $U$ is defined on $\widehat{M} = M\backslash\Delta$, with 
collision set 
\[\Delta = \{q\in M\mid q_i=q_j\mbox{ for some }i\neq j\}.\]
In Section \ref{section:notation} the stage is set by introducing the notation necessary to
treat the combinatorics encoded in $\Delta$.

The $n$ particle phase space $\widehat P = T^* \widehat{M}$ over noncollision
configuration space $\widehat{M}$
is non-compact in three ways: 
\begin{enumerate}[$\bullet$]
\item 
Total energy $H: \widehat P\to \bR$ may have any value $E\in \bR$.
\item 
There are collisions between the particles, with limits in $\Delta$.
\item
Particles may escape to spatial infinity. 
\end{enumerate}
To  compactify phase space 
the first obvious step\footnote{We don't touch upon the question about the limits $E\to\pm\infty$.} 
is to separately consider motion on the 
energy surfaces $\hSE$.
In Section \ref{sub:blow:conf} 
we compactify configuration space $\widehat{M}$ by adding a sphere $\bS$ at infinity to $M$
and by forming the real blow up of the  collision set $\Delta$,  
thus obtaining a {\it compact}  manifold with corners $\wideparen{M}$.
In Subsection \ref{subsec: rescale vel} we rescale velocities (or momenta)  in order to 
compactify the energy surface $\wideparen{\Sigma}_E$.  
In Subsection \ref{subsec: rescale time} we
rescale time in order to get a well-defined dynamics on  the compactified energy surface.  
We will see that $\wideparen{\Sigma}_E = \wideparen{M}_E \times \bS$
where $\wideparen{M}_E \subseteq \wideparen{M}$ is the compactified Hill region 
$\{q: U(q) \ge - E \}$.
In particular, points along  the traditional Hill boundary have been blown up into spheres.
As a result, the original dynamics has been altered at the Hill boundary
and some care is required to understand  the rescaled dynamics there and recover the original dynamics.  
By continuity the flow on the boundary mimics the one near the
boundary, providing useful information about the original flow 
(see Remark \ref{dynamic philosophy}).\\
Finally, in Section \ref{sec:topology} the topology of the construction is considered.

\subsection{History and Comparisons}
An incomplete flow on a noncompact manifold $X$ can always be compactified, but in a  
useless   way.  Choose an  arbitrary compactification $\bar X$ of $X$.  
(There are  many! \footnote{
 any compact connected $N$-manifold arises as the  compactification
of Euclidean $N$-space})
Slow down the flow
by rescaling the vector field which defines this flow so that the vector field vanishes on 
$\bar X \setminus X$.  Voila!  (We thank P. Deligne for this remark.)
The real trick is to compactify $X$ and rescale the vector field  in such a way 
that the extended rescaled field  {\bf does not vanish} on 
$\partial X$ and  allows  you to extract new information
about the original flow on $X$.

McGehee's partial compactification    \cite{McG1} succeeds eminently for extracting
information  about near-total collision motions.   
{\sc Robinson} \cite{Ro}, {\sc McGehee} \cite{McG3} and others have   
``turned McGehee's microscope
around to become a  telescope''  by adding a 
a manifold  at infinity, instead of at collision, this infinity manifold being used  to account for escape orbits.  
  Again, they derived non-trivial information
 about the dynamics by this trick. {\bf In essence, what we do here is systematically
 implement both of these partial compactifications in order to get our full compactification.}  
 
Ours is the first paper to fully compactify  classical  $n$-body problems, $n> 2$, 
in a {\it potentially useful way}.   (We have yet to establish its utility.) 

Two-body problems have been compactified. 
Moser in particular compactifies
the $\alpha =1$ negative energy 2-body problem in $d$ dimensions, 
in the center of mass frame. 
Moser's  compactified phase space is the unit sphere bundle of the $d$-sphere. 
Like the methods of {\sc Levi-Civita} \cite{LC},  
Kuustanheimo-Steifel and a number of others, Moser's method
is a `regularization' rather than a `blow-up' and as such is quite particular to  
$\alpha =1$ and so we will reject it.   
 
 Extended versions  of regularizations, specific for $\alpha =1$,   were developed by 
{\sc Heggie} \cite{He},  {\sc Waldvogel} \cite{Wald},
{\sc Lema\^{i}tre}, \cite{Le},  {\sc Moeckel-Mont\-gomery} \cite{MM}, and others
for improving  numerics and deriving partial compactifaction results for 
 $n$-body problems with $\alpha =1$.

   {\sc Qiu-Dong Wang}  \cite{Wa} uses a time
and velocity rescaling very similar to ours
  but for different purposes.    In particular, he does not add in collision manifolds
or manifolds at infinity.  By adding in these manifolds
we allow for coherent $\alpha$ and $\omega$ limit sets for collections
of  orbits leaving the honest energy $E$ phase space $\wideparen{\Sigma}_E$.
These limit sets are actual 
places for orbits to go to and will afford us,   we hope, with an eventual better understanding
of the near-collision orbits and of the way in which   clusters of particles approach spatial infinity and
to what extend the clusters asymptotically become independent. 
See the  Remark \ref{dynamic philosophy}.

{\sc Graf} \cite{Gr} constructed what we nowadays call {\em Graf partitions}
in order to prove asymptotic completeness for the quantum $n$-body problem.
 {\sc Vasy}, a student of Melrose, used in  \cite{Va} essentially these same 
Graf partitions  combined with 
Melrose style blow-ups into manifolds-with-corners in order to obtain 
new information about scattering in quantum 3-body and $n$-body problems.

\subsection{Goals}  The flows on the boundary strata are simpler 
than the flow on the bulk  ($\widehat{\Sigma}_E$).
For example, the flow on the open part of the locus at infinity is a reparameterization of 
geodesic flow on the sphere, while the flows on the open strata of the 
binary collision locus are reparameterized Kepler flows.   

Let us agree that we can concatenate two boundary trajectories if the 
$\omega$-limit of one agrees with the $\alpha$-limit of the other.
Thus, if $\gamma_1$ and $\gamma_2$ are boundary trajectories for which 
$\lim_{t \to + \infty} \gamma_1 (t) = \lim_{t \to -\infty} \gamma_2 (t) $
then we can form the concatenation $\gamma_2 * \gamma_1$ 
made up of first traversing $\gamma_1$
and then traversing $\gamma_2$.     

We believe that these concatenated boundary trajectories form a kind
of skeleton, or support locus which controls  the flow in the bulk near the boundary.  
Specifically,  given any such
concatenation $c = \gamma_2 * \gamma_1$ we believe that we can
prove the existence of  sequences $c_i$ of trajectories which lie completely 
in the bulk and which converge to $c$
in an appropriate sense: $\lim_{i \to \infty} c_i = c$.
Establishing the existence of these ``shadowing sequences'' 
 $c_i$  is work in progress.  It would
validate various observations and numerical experiments made in  \cite{DMMY} and \cite{FKM}.

\bigskip
\noindent
{\bf Acknowledgement:}\\ 
We dedicate this article to Alain Chenciner, our guiding star, at his 80th birthday.\\[2mm]
AK thanks Eva Miranda (Barcelona) for motivating discussions.\\[1mm]
RM thanks the Simons foundation for travel support. \\[1mm]
We both thank the anonymous referee for useful suggestions.
%
\section{Notation} \label{section:notation}
%
\subsection{Phase flow}
%
Consider $n$\label{number of particles} 
particles of masses $m_i>0$\label{masses}  moving in $d$-dimensional Euclidean space $\bR^d$.
(Take $n > 1$ please!) 
Introduce a separate copy 
$M_i:=\bR^d$ of the Euclidean  space for each particle $\ i\in N := \{1,\ldots,n\}$ \label{set of particles}
\label{configuration space of i--th particle} so that  $q_i \in M_i$.  The center of mass zero configuration space  \label{total configuration space} is
\[\textstyle
M := \big\{q\in\bigoplus _{i\in N} M_i \mid \sum_{i\in N} m_i q_i = 0\big\}  \]
and forms a codimension $d$  linear subspace of the vector space $(\bR^d)^n$.
We   use  the mass-inner product $\LA q,q' \RA_\cM = \sum_i m_i q_i \cdot q_i '$ on $M$
instead of  the standard   inner product  
$\LA q , q' \RA = \sum_i q_i \cdot q_i ' $.   
The  mass matrix $\cM: (\bR^d)^n \to (\bR^d)^n$ defined by 
 $\cM(q_1, \ldots , q_n) = (m_1 q_1, m_2 q_2, \ldots , m_n q_n)$
intertwines the two inner products:  
\beq
\LA q,q'\RA_\cM:=  \LA q,\cM q'\RA \, . 
\eeq
We set $\|q\|_\cM:=\sqrt{\LA q,q\RA_\cM}$ 
omitting the subscript $\cM$ whenever possible.  We can write
the kinetic energy as 
\beq 
K(p) := \textstyle{\frac{1}{2}} \LA p, \cM^{-1} p \RA \, .
\eeq
The {\em collision set} in configuration space is given by
\beq
\Delta := \{q\in M\mid q_i=q_j\mbox{ for some }i\neq j\in N\} \, .
\Leq{coll:set}
We consider homogeneous potentials 
on the {\em noncollision configuration space}\label{noncollision configuration space}\beq
\widehat{M} := M\backslash\Delta .
\Leq{hat:M}
On the phase space $\hP:=T^*\widehat{M}$ 
the Hamiltonian function is given by
\beq
H:\hP \to \bR \qmbox{,} H(q,p) = K(p)- U(q) \, ,
\Leq{Y}
with the real-analytic potential $U: \widehat{M}\to \bR$ of the form \eqref{Ham:Pot},
with $Z_{i,j} > 0$.

With the Euclidean gradient $\nabla_{\!(i)}$ on each $M_i = \bR^d$ 
the Hamiltonian equations of \eqref{Y} have the form 
\beq
\dot{q}_i=\frac{p_i}{m_i}\qmbox{,}
\dot{p}_i= \sum_{j\in N\setminus \{i\}} \nabla_{\!(i)} U_{i,j}(q_j-q_i)
\qquad (i\in N).
\Leq{Ham:eq}
where
\beq  
U_{i,j}(q_j-q_i) = \frac{Z_{i,j}}{\| q_i - q_j \|^{\alpha}}
\Leq{two-particle potl: eq}
A useful alternative way to write the equations is in terms of velocities: 
\beq
\dot q = v \qmbox{ , }   \dot v = \nabla U(q)
\Leq{N:eq}
where now $q, v \in M$, $q \notin \Delta$ 
and the gradient $\nabla$ is with respect to the mass metric.
\begin{remark}[other potentials]\quad\\ 
The same form for Newton's equations holds  for any potential  
$U (q_1, \ldots , q_n) = \sum_{i<j} U_{i,j} (q_i - q_j)$
which is a sum of pair potentials $U_{i,j} : \bR^d \setminus \{0 \} \to \bR$.\\
Many of our results leading up to our theorem will also hold for these more general potentials provided their pair potentials have
appropriate blow-up and decay conditions mimicking that of the power law potentials. 
One notable such required condition would be  
$U_{i,j} (z) \sim Z_{i,j} /\|z\|^{\alpha} + O(\|z\|^{1- \alpha})$  with $0 < \alpha < 2$ 
along with corresponding conditions on the  derivatives of $U_{i,j}$
as $z \to 0$.  It will be important that $\alpha$ does not depend on $i,j$.  
\hfill $\Diamond$
\end{remark} 
Using the natural symplectic form $\omega_0$ on $\hP=
T^*\widehat{M}$ we write \eqref{Ham:eq} in the form
$\dot{x}=X_H(x)$ for the Hamiltonian vector field
$X_H$ defined by ${\bf i}_{X_H}\omega_0=dH$.

The flow of this vector field is real-analytic and fixes energy so
defines a flow on each of the {\em energy surfaces}\label{energy surface}
\beq
\hSE := \big\{x\in\hP\mid H(x)=E \big\}\qquad(E\in\bR).
\Leq{energy:surface} 
%
\subsection{Cluster decomposition}
%
Cluster decompositions provide us with  the  book-keeping we need to
index the  ways we can end  on the collision locus
$\Delta$.  We borrow the language from   combinatorial theory.  See {\sc Aigner} \cite{Ai}.
\begin{defi}\quad\label{def:meet:join}\\[-6mm]
\begin{enumerate}[$\bullet$]
\item
A {\bf partition} or {\bf cluster decomposition}
of $N = \{1, 2, \ldots , n \}$ is a collection 
$\cC:=\{C_1,\ldots,C_k\}$ of disjoint non-empty subsets of $N$ whose union is $N$.
The elements of $\cC$ are called the {\bf atoms} or {\bf clusters} of the cluster decomposition. 
\item
A cluster decomposition $\cC$ induces an  equivalence relation on $N$ whose equivalence classes are $\cC$'s atoms.
We write $[i]_\cC$ or simply  $[i] \in \cC$ for the atom containing $i \in N$.
\item
The {\bf partition lattice} $\cP(N)$ is the set of 
cluster decompositions $\cC$ of $N$, partially ordered by 
{\bf refinement}, that is
\[\cC=\{C_1,\ldots,C_k\}\preccurlyeq \{D_1,\ldots,D_\ell\} = \cD, 
\qmbox{if} C_m\subseteq D_{\pi(m)}\] 
for a suitable surjective relabelling  map
$\pi:\{1,\ldots,k\}\to\{1,\ldots,\ell\}$. Then $\cC$ is called {\bf finer} than $\cD$ and $\cD$ 
{\bf coarser} than $\cC$.
\item
The {\bf rank} of $\cC\in\cP(N)$ is the number $|\cC|$ of its atoms.
\item
The   finest and coarsest  elements of $\cP(N)$ are denoted by
\beq
\cC_{\min} := \big\{ \{1\},\ldots,\{n\} \big\}\qmbox{and}
\cC_{\max} := \big\{ \{1,\ldots,n\}\big\} \, ,
\Leq{fine:coarse}
and we set
\beq
\cP_\Delta(N) := \cP(N)\setminus\{\cC_{\min}\}
\qmbox{and}
\cP_\bS(N) := \cP(N)\setminus\{\cC_{\max}\} \, .
\Leq{cP:Delta}
\end{enumerate}
\end{defi}
For a partition $\cC$ we define the $\cC$--{\em collision} subspace
\[\De_\cC := \l\{q\in M\mid q_i=q_j\ \mbox{if}\ [i]_\cC=[j]_\cC \ri\} \, .\]
Note that $\De_\cC \subseteq \Delta$
as long as $\cC \ne \cC_{\min}$.  
We have that
\[ \cC \preccurlyeq   \cD \ \implies \ \De_\cD \subseteq \De_\cC\]
and that 
$$\De_\cC=\bigcap_{C\in\cC}\De_{C}$$
where, for a subset $C\subseteq N$ we declare  
\beq
\De_C := \{q\in M\mid q_i=q_j\ \mbox{for}\ i,j\in C\} \, .
\Leq{def:MEC}
The superscript $E$ refers to the fact that the `non-frozen'   coordinates are 
{\em external} to those associated to the atom $C$ of $\cC$.

For subsets $C \subseteq N$ we denote  the 
$\cM$--{\em orthogonal} complement to $\De_C$ by 
\beq  
\Di_C : = (\De_C)^{\perp} \,.
\Leq{def: MiC}
One computes
\beq \textstyle
\Di_C = {\rm Ker}\l(\Pi^E_C\ri)
=\big\{q\in M\mid q_i=0\ \mbox{for}\ i\not\in C,\ \sum_{i\in C}\,m_iq_i =0\big\} \, .
\Leq{M:c:bot}
In the first equality of \eqref{M:c:bot},  
$\Pi^E_C$\label{cluster projection} denotes the orthogonal 
projection onto $\De_C$. The superscript $I$ refers to the fact that only coordinates
{\em internal} to the atom $C$ can vary on this subspace.  Thus the orthogonal projection
onto $\Di_C$ is  $\Pi_C^I:=\idty_M-\Pi^E_C$.

Similarly for a partition $\cC$ we define
$$
\Di_\cC : = (\De_\cC)^{\perp} \, .
$$
Then $\Di_\cC = \bigoplus_{C\in \cC} \Di_{C}$,
since $\De_\cC=\bigcap_{C\in \cC} \De_{C}$. By (\ref{M:c:bot})
\beq
M = \De_\cC \oplus \Di_\cC = \De_\cC \oplus \bigoplus_{C\in \cC}\Di_{C}
\Leq{decomp}
is an $\cM$--orthogonal decomposition. 
Associated to the orthogonal decomposition we have the  orthogonal projections  
\beq
\Pi^E_\cC := \prod_{C\in\cC}\Pi^E_C
\qmbox{, respectively} \Pi^I_\cC := \idty_M-\Pi^E_\cC=\sum_{C\in\cC}\Pi^I_C \, .
\Leq{proj:E} 
One easily computes the dimensions
\beqn
\dim(\De_\cC) &=& \textstyle
d\big(n-1-\sum_{C\in \cC}\ (|C|-1)\big) =  d\,(|\cC|-1)\, ,\nonumber\\
\dim(\Di_\cC) &=& \sum_{C\in \cC}\dim(\Di_{C}) 
= d\sum_{C\in \cC}(|C|-1) 
= d(n-|\cC|) \, .
\label{dim:MC}
\eeqn
To lighten the notation set
\[ q^E_\cC := \Pi^E_\cC(q)\qmbox{and}q^I_\cC:=\Pi^I_\cC(q)\qquad(q\in M).\]
We will  even omit the subscript $\cC$ when the context permits.
For a nonempty subset $C\subseteq N$ we define the \label{cluster mass} 
{\em cluster mass\,} and {\em cluster barycenter} of $C$ by
\[m_C := \sum_{j\in C} m_j \qmbox{and}
  q_C := \frac{1}{m_C} \sum_{j\in C} m_j q_j \, .\]
In particular $m_N$ equals the 
{\em total mass} of the particle system.
Then for the partitions $\cC\in \cP(N)$ the $i$--th component of the 
cluster projection is given by the barycenter
\beq
\l(q^E_\cC\ri)_i = q_{[i]_\cC} \qquad(i\in N)
\Leq{cl:bar}
of its atom. Similarly\label{Deltaqi}
\[ \l(q^I_\cC\ri)_i = q_i-q_{[i]_\cC} \qquad (i\in N)\]
is its distance from the barycenter. 
The {\em scalar moment of inertia}\label{moment of inertia}
\beq
J:M\to\bR \qmbox{,} J(q):=\LA q,q\RA_\cM
\Leq{scalar:moment:of:inertia}
splits into the {\em cluster barycenter moment}
\[J^E_\cC:=J\circ\Pi^E_\cC \qmbox{,} J^E_\cC(q)=\sum_{C\in\cC} m_{C}\LA q_C,q_C\RA\]
and the {\em relative moments of inertia} of the clusters $C\in\cC$
\[J_{C}^I := J\circ\Pi_{C}^I \mbox{ , } 
J_{C}^I(q) = \sum_{i\in C} m_i\|\l(q^I_\cC\ri)_i\|^2 = 
\frac{1}{2m_C} \sum_{i,j\in C} m_im_j\LA q_i-q_j,q_i-q_j\RA ,\]
that is
\beq
J=J^E_\cC+J^I_\cC\qmbox{for} 
J^I_\cC(q) := \sum_{C\in\cC}J_{C}^I(q) = \LA q^I_\cC,q^I_\cC\RA_\cM.
\Leq{J:decomp}

\begin{example}[a binary pair cluster and Jacobi vectors]\label{ex:1}\quad\\
For a  pair $i\ne j$ of particle labels we can form
the cluster decomposition $\cC$ of rank $n-1$ whose only non-singleton atom is
the cluster $C:=\{i,j\}$.  These partitions correspond to isolated binary collisions
and   generate $\cP(N)$ under join.  
For simplicity of notation we will  take $n =3$ and   $C = \{ 1,2 \}$. 
$M$ has dimension $2d$,  $\Delta^E _{(3)} = M$, $\Delta^I_{(3)} = 0$ while 
$\Delta^E _{(1,2)}$ and $\Di_{(1,2)}$ are both $d$--dimensional.  
We have $\De_\cC = \De_{(1,2)}, \Di_\cC = \Di_{(1,2)}$ 
and  $\De_{(1,2)} \oplus \Di_{(1,2)} = M$ (orthogonal direct sum). 
We compute  
\[(\qe_\cC)_1 = \frac{m_1 q_1 +m_2 q_2}{m_1 + m_2}  = (\qe_\cC)_2 \qmbox{,}  
(\qe_\cC)_3 = q_3 \] 
while 
\[(\qi_\cC)_1 = \frac{m_2}{m_1+m_2}(q_1-q_2) \qmbox{,} (\qi_\cC)_2 = \frac{m_1}{m_1+m_2}(q_2-q_1)\qmbox{,} (\qi_\cC)_3 = 0 \, .\]
The vector $\qi_\cC$ is parameterized by the   
Jacobi vector  $\xi_1 = q_1 - q_2 $. The vector $\qe_\cC$ is parameterized by   the  other 
Jacobi vector   $\xi_2 =  q_3 - \frac{m_1 q_1 +m_2 q_2}{m_1 + m_2}$.  (Use $m_1 q_1 + m_2 q_2 + m_3 q_3  = 0$
to show that $\xi_2, q_3$ and  $\frac{m_1 q_1 +m_2 q_2}{m_1 + m_2}$
are all non-zero scalar multiples of each other.)   The identity (\ref{J:decomp})
becomes the traditional quadratic decomposition   
\beq
\| q \|^2 = \mu_1 \| \xi_1 \|^2 +  \mu_2 \|\xi_2 \|^2
\Leq{J:K:I:mass}
with mass coefficients $\mu_i$ given by 
\[ \textstyle
\frac{1}{ \mu_1}= \frac{1}{m_1} + \frac{1}{m_2}  \qmbox{,} 
\frac{1}{ \mu_2} =  \frac{1}{m_1 +m_2} +  \frac{1}{m_3} \, .
\hfill\tag*{$\Diamond$}\]
\end{example}
%
\subsection{On to phase space} 
%
Let $M^*$ denote the dual space of our vector space $M$.  There  are natural 
identifications $TM\cong M\times M,\ T^*M\cong M^*\times M$ of the tangent
space resp.\ phase space of $M$. These gives rise to the inner products
\[\LA \cdot,\cdot\RA_{TM}:TM\times TM\to\bR \qmbox{,} 
\LA (q,v),(q',v')\RA_{TM}:=\LA q,q'\RA_\cM+\LA v,v'\RA_\cM\]
and
\beq
\LA \cdot,\cdot\RA_{T^*M}:T^*M\times T^*M\to\bR \qmbox{,} 
\LA (q,p),(q',p')\RA_{T^*M}:= \LA q,q'\RA_\cM+\LA p,p'\RA_{\cM^{-1}}
\Leq{cotangent:scalar:product}
$\l(\mbox{with}\ \LA p,p'\RA_{\cM^{-1}}=\sum_{i=1}^n\frac{\LA p_i,p_i'\RA}{m_i}\mbox{ for 
the {\em momentum vector} } p=(p_1,\ldots,p_n)\ri)$.

The tangent space $TU$ of any linear subspace $U\subseteq M$ is naturally a 
linear subspace of $TM$. 
Using the inner product, we can also consider $T^*U$ as a subspace of $T^*M$.

We thus obtain $T^*M$--orthogonal decompositions
\[T^*M =
T^* (\De_\cC) \oplus\bigoplus_{C\in\cC}T^*(\Di_C)\qquad
\big(\cC\in\cP(N)\big)\]
of phase space. 
With  
\[\widehat{\Pi}^I_\cC := \idty_{T^*M}-\widehat{\Pi}^E_\cC
=\sum_{C\in\cC}\widehat{\Pi}^I_C\]
the $T^*M$--orthogonal projections 
$\widehat{\Pi}^E_\cC, \, \widehat{\Pi}^I_\cC: T^*M \to T^*M$
\label{widehatPicC}
onto these subspaces are given by the {\em cluster coordinates}
\label{cluster coordinates}
\beq
(\qe,\pe) := \widehat{\Pi}^E_\cC(q,p) \qmbox{with} (\qe_i,\pe_i)
= \l(q_{[i]} , \frac{m_i}{m_{[i]}}p_{[i]} \ri)\quad(i\in N),
\Leq{cl:co}
and {\em relative coordinates}\label{relative coordinates}
\[(\qi, \ppi) := \widehat{\Pi}^I_\cC(q,p) \qmbox{with} 
(\qi,\ppi_i) = (q_i-\qe_i , p_i-\pe_i ) \quad (i\in N).\]
Here $p_C := \sum_{i\in C} p_i \in \bR^d$ is the {\em total momentum of the cluster}
$C\in \cC$. Unlike in (\ref{cl:bar}) we omitted the subindex $\cC$
in (\ref{cl:co}), but will include it when necessary.

With this notation the equations of motion for particle no.\ 
$i\in N$ are
\[\frac{d}{dt} \qe_i = m_i^{-1} \pe_i\qmbox{,}
\frac{d}{dt} \pe_i = \frac{m_i}{m_C} \sum_{j,k\in
N:[j]=[i],\ [k]\neq[i]}\nabla U_{j,k}(q_j-q_k) \, ,\]
$\frac{d}{dt} \qi_i = m_i^{-1} \ppi_i$ and
\beq
\frac{d}{dt} \ppi_i = \sum_{k\in N\setminus{\{i\}}} \nabla U_{i,k}(q_i-q_k)
-\frac{m_i}{m_C} \sum_{j,k\in
N:[j]=[i],\ [k]\neq[i]}\nabla U_{j,k}(q_j-q_k) \, . 
\Leq{dqpI:dt}
\begin{lemma}\label{lem:sym}
The vector space automorphisms 
\beq
\l(\widehat{\Pi}^E_\cC,\widehat{\Pi}^I_\cC \ri):T^*M\longrightarrow
T^*\De_\cC\oplus\bigoplus_{C\in\cC}T^*(\Di_C)
\qquad\big(\cC\in \cP(N)\big)
\Leq{sy:trans}
are symplectic w.r.t.\ the natural symplectic forms on these cotangent bundles.
\end{lemma}
{\bf Proof.}
This follows from 
$T^*\big(\De_\cC\oplus\bigoplus_{C\in\cC}\Di_C\big)=T^*\De_\cC\oplus\bigoplus_{C\in\cC}T^*(\Di_C)$.\hfill$\Box$\\[2mm]
The {\em total kinetic energy}\label{total kinetic energy}
\[ \textstyle
K:T^*M\to\bR \qmbox{,} K(q,p) \equiv K(p)=\eh\LA p,p\RA_{M^*} =
\sum_{i=1}^n\frac{\LA p_i,p_i\RA}{2m_i}\]
 splits into the  external or {\em barycentric} kinetic energy \label{barycentric kinetic energy}\label{KcC} 
\beqno
K^E_\cC &:=& \textstyle
K\circ\widehat{\Pi}^E_\cC \qmbox{,} K^E_\cC(q,p) = \sum_{C\in\cC}
\frac{\LA p_C,p_C\RA}{2m_C}
\eeqno
and internal or {\em relative} kinetic energy associated to each   cluster $C\in\cC$
\beqno
K_C^I &:=& \textstyle
K\circ\widehat{\Pi}_C^I \qmbox{,} 
 K_C^I(q,p) = \sum_{i\in C}\frac{\LA \ppi_i,\ppi_i\RA}{2m_i} \, .
\eeqno
That is,
\[  K = K^E_\cC + K^I_\cC \qmbox{with}K^I_\cC:= \sum_{C\in\cC}K_C^I \, .\]
The {\em internal} and  {\em external cluster potentials} 
and  {\em Hamiltonians} are given by
\beq
U^I_\cC(q) := \sum_{C\in \cC} U^I_C(q)
\qmbox{with}
U^I_C(q) := \sum_{i<j\in C} U_{i,j}(q_i-q_j),
\Leq{VcC:I}
\beq
H^I_\cC(p,q) := K^I_\cC(p) - U^I_\cC(q) = \sum_{C\in \cC} H^I_C(p,q)
\qmbox{with}H^I_C(p,q) := K^I_C(p) - U^I_C(q)
\Leq{HcC:I}
and
\beq
U^E_\cC(q) :=
\sum_{i<j\in N:\; [i]_\cC \neq [j]_\cC} U_{i,j}(q_i-q_j)
\qmbox{,}H^E_\cC(p,q):= K^E_\cC(p) - U^E_\cC(q).
\Leq{VcC:E}
We have for all $\cC\in\cP(N)$
\beq
U = U^I_\cC+U^E_\cC \qmbox{and} H = H^I_\cC+H^E_\cC.
\Leq{UUU}
(Unlike the kinetic energies,   the cluster potentials cannot be written as
$U\circ \Pi^I_\cC$ etc).
\begin{remark}[partition of configuration space]\quad\label{strata}\\
The linear subspaces $\De_\cC$  
generate a   stratification\,\footnote{A stratification of a manifold is a locally finite
partition into smooth submanifolds.  See for example \cite[p. 83]{Wall}. }
of $M$ with strata  
\beq
\Xi_\cC \equiv\, \Xi_\cC^{(0)} := \De_\cC\Big\backslash
\bigcup_{\cD\succneqq\cC}\De_\cD\qquad\big(\cC\in \cP(N)\big).
\Leq{Xi:Null}
The upper index $(0)$ is generalized in \eqref{Z4} but omitted if there is no danger of 
confusion.
$\Xi_\cC$ consists of those  collisions where only 
those  particles whose particle indices belonging to the same $C_\ell \subseteq N$
of   $\cC = \{C_1,\ldots,C_k\}$ coincide.
As $\Xi_\cC$ is open in the vector space $\De_\cC$, by \eqref{dim:MC} it is a manifold with 
dimension $\dim(\Xi_\cC) = d\,(|\cC|-1)$.
\hfill $\Diamond$
\end{remark}
%
\section{Blowing up the configuration space}\label{sub:blow:conf}
%
The center of mass configuration space $M$, being a  finite-dimensional  real vector space,  is diffeomorphic to its
open unit ball $B\subseteq M$.   For diffeomorphism we can take  
\[ \Phi: M\to B\qmbox{,} 
\Phi(0)=0\ ,\ \textstyle 
\Phi(q) = \tanh(\|q\|) \frac{q}{\|q\|} \quad (q\in M\backslash \{0\} ) \,.\] 
$B$ is compactified by attaching its
boundary $\bS:=\partial B$, the unit sphere of dimension $d(n-1)-1$ 
which corresponds to letting
$\|q \| \to \infty$ in the expression for $\Phi$.  
Abusing notation, we set $\partial M:=\bS$, and 
\beq
\overline{M} := M\sqcup \bS \, ,
\Leq{closed:ball} 
a manifold with boundary.
In this sense the topological boundary of the
open configuration manifold $\widehat{M} = M\backslash\Delta$ is $\Delta\sqcup \bS$.

The partial compactifications 
$\wideparen{M}_{\widehat \bS}$ of $\widehat{M}\sqcup \widehat \bS$ along 
$\widehat{\bS} := \bS\backslash \Delta$ (Subsection \ref{sub:conf:infty}) and
$\wideparen{M}_\Delta$ of $\widehat M$ along $\Delta$ 
(Subsection \ref{sub:conf:delta}) 
are of independent interest. 
The first  is relevant for the dynamics of the $n$ 
particles as their  mutual distances go
to infinity and  was treated in \cite{DMMY}.
The second corresponds to collisions and  leads to a generalisation of the blow-up of the
total collisions for the case $n=3$ as first
treated by {\sc McGehee} in \cite{McG1}.
 
The full compactification $\wideparen{M}$ in Subsection \ref{sub:conf:compact} 
gives rise to an additional aspect, because
it includes the points of the sphere $\bS$ at infinity corresponding to non-trivial
clusters.
This will, we hope (work in progress)   lead us  to a positive solution of the problem of 
asymptotic completeness\footnote{See 
{\sc  Derezi\'{n}ski} and {\sc G\'erard} \cite[Section 5.10]{DG}.}
in the 
three-body problem, thus asymptotically relating the joint motion of the particles and the 
pure two-body dynamics in the non-trivial cluster.\footnote{In the case of smooth 
{\em bounded} pair 
potentials that are $\alpha$-homogeneous for large distances, a positive solution of the
asymptotic completeness problem for arbitrary number $n$ of particles is expected, too.
In the present setting, however, the motion is not asymptotically complete 
for $n>3$ because of the existence of non-collision singularities.} 
%
\subsection{Real blow up,  generally}\label{sub:conf:gen}
%
We will be repeatedly implementing the real blow-up construction
as described in Melrose \cite[Sect.\ 5.3]{Me}.  Our constructions will be concrete
and essentially linear-algebraic so the general construction is not needed here.
However, it may be useful to get a rough understanding of it, if for
no other reason than to familiarize ourselves with the notation.

The general construction proceeds as follows.  
Given an embedded submanifold $Y$ of a manifold $X$
the real blow-up $[X:Y]$ is formed by deleting $Y$ from $X$ 
and replacing it with the space $S^+N Y$ of rays in the  normal bundle $NY$ to $Y$.
This space $S^+N Y$   is a sphere bundle over $Y$ 
with the spheres having one less than the codimension of $Y$.
The resulting $[X:Y]$  is a manifold with boundary, that boundary
being the sphere bundle.  
The manifold  comes with a smooth blow-down map $[X:Y] \to Y$ which
takes the sphere bundle onto $Y$ by the bundle projection 
and is a diffeomorphism away from the sphere bundle.
The manifold structure is obtained by invoking the tubular neighborhood theorem.  
This is commonly done by  using an auxiliary Riemannian structure which allows
us to use geodesics normal to $Y$ to  form a diffeomorphism between a 
neighborhood of $Y$ and a neighborhood of the
zero section of the normal bundle of $Y$.  
When one looks at things in local Gaussian-cylindrical coordinates
about $Y$ the whole construction boils down to using polar coordinates normal to $Y$.

The construction works if $X$ has a boundary and $Y$ is transverse to the boundary, 
intersecting it in its own boundary.
In that case $[X:Y]$ is a manifold with codimension two corners.  
The construction can be iterated upon choosing a finite collection
$Y_a, a \in I$ of embedded submanifolds, provided certain conditions are verified 
concerning the intersections of the
closures of the $Y_a$  with each other, and with corners arising in previous steps,  
yielding manifolds with  deeper and deeper corners.  
\subsection{Blowing up configuration space at infinity}\label{sub:conf:infty}
%
Although we do not blow up the entire configuration vector space $M$, 
diffeomorphic to the open ball $B\subseteq \bR^k$ of radius one, it is instructive
to notice that this would just  reproduce the closed ball $\overline{B} = B\sqcup \bS$
from \eqref{closed:ball}.  In other words, 
\[ [\,\overline{M}: \bS\,] \ \cong \  \overline{M} \,.\]
The reason is the following one:
\begin{example}[blowing up configuration space $M$]\quad\\  
As $\bS$ is the boundary of the configuration manifold $M$, its blow up is based on 
the general definition given in {\sc Melrose} \cite[Sect.\ 5.3]{Me}.
As $\overline{M}\backslash \bS = M$, we thus set
\[ [\overline{M}: \bS] := M \sqcup (S^+\!N \,\bS) \, ,  \]
with $S^+N\, \bS$ being the {\em inward pointing part} of the normal sphere bundle of 
$\bS\subseteq \overline{M}$. Since $\bS\subseteq \overline{M}$ is of codimension one,
$S^+\!N \,\bS$ is diffeomorphic to $\bS$, so that we get a simple result of that blow up:
$[\overline{M}: \bS]\cong  \overline{M}$, the closed unit ball.
\hfill $\Diamond$
\end{example}
With the open subset $\widehat{\bS} := \bS\backslash \Delta$ of the sphere 
we similarly obtain the disjoint union 
\beq 
\wideparen{M}_{\widehat \bS} \
:= [\widehat{M}\sqcup \widehat{\bS}:\widehat{\bS}]  \cong \widehat{M}\sqcup \widehat{\bS}\, ,
\Leq{wideparen:M:S}
which is a manifold with boundary.
\subsection{Blowing up configuration space at collisions}\label{sub:conf:delta}
%
Next we blow up the configuration space $\widehat{M} = M\backslash\Delta$ along 
the thick diagonal
\[\Delta=\bigcup_{\cC\in \cP_\Delta(N)}\  \De_\cC \, ,\]
using the family \eqref{decomp} of $\cM$-orthogonal decompositions 
$M = \De_\cC \oplus \Di_\cC$ (Recall  that $\cP_\Delta(N)$ denotes 
$\cP(N) \setminus \{\cC_{\min}\}$.)
With the definition \eqref{J:decomp} of $J^I_\cC$, we write the coordinate 
$q^I_\cC$ in the form
\beq
q^I_\cC = r \,Q^I_\cC\quad\mbox{with } J^I_\cC(Q^I_\cC)=1 \mbox{ and } 
r:=\big(J^I_\cC(q^I_\cC) \big)^{1/2}\qquad \big(\cC\in \cP_\Delta(N)\big).
\Leq{q:r:Q}
For $q^I_\cC\neq 0$ this polar decomposition is unique and
\beq
\widehat{\Delta}^I_\cC := \Delta^I_\cC\setminus\{0\} \cong (0,\infty)\times S^I_\cC \quad 
\mbox{, with the sphere }
S^I_\cC:=\big(J^I_\cC\big)^{-1}(1) \,.
\Leq{def:C:sphere}
Somewhat loosely speaking, we call the  $Q^I_\cC$ of  equation (\ref{q:r:Q}) the  
{\em coordinates} on~$S^I_\cC$.\\
We consider $M\cong \bR^{(n-1)d}$ as the vector bundle 
$M\cong \De_\cC \oplus \Di_\cC\to \De_\cC$
and follow {\sc Melrose} \cite[Sect.\ 5.2]{Me}  
in defining the blow-up of $M$ along the zero section 
$\De_\cC \times \{0\} \cong \De_\cC$
as the manifold with boundary
\[ [M:\De_\cC] := \De_\cC\times \wideparen{\Delta}^I_\cC \qmbox{with}
\wideparen{\Delta}^I_\cC:=[0,\infty)\times S^I_\cC \, .\]
The diagonal map
\beq
\hat{I}_\Delta:\widehat{M} \longrightarrow  
\prod_{\cC\in \cP_\Delta(N)} [M:\De_\cC] \qmbox{,}q\longmapsto (q)_{\cC\in \cP_\Delta(N)}
\Leq{hat:I:Delta}
smoothly imbeds $\widehat{M}$ as a submanifold of 
the manifold $\prod_{\cC\in \cP_\Delta(N)}  \De_\cC\times \widehat{\Delta}^I_\cC$.

\begin{figure}[h]
\centerline{\includegraphics[width=100mm]{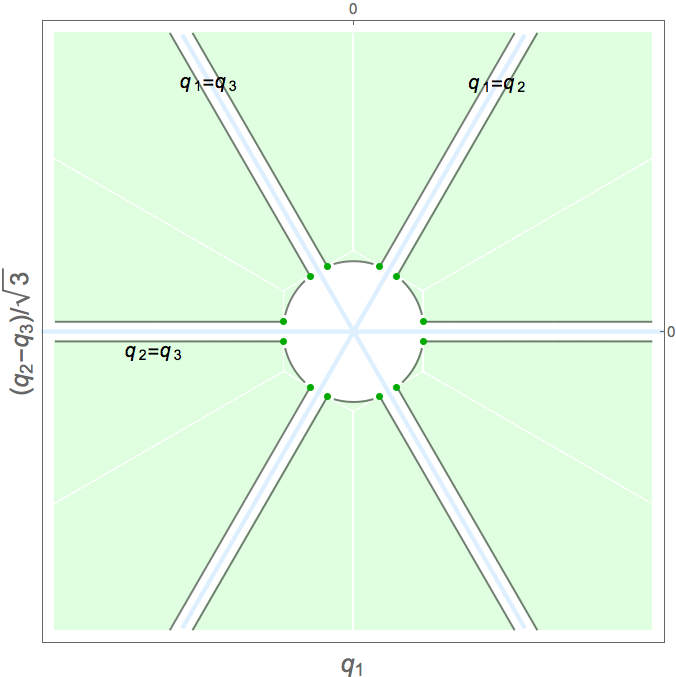}}
\caption{Schematic view of the manifold with corners $\wideparen{M}_\Delta$  (light green)
arising by blowing up configuration space 
$M \cong \bR^{(n-1)d}$, for $n=3$ particles in $d=1$ dimension. 
The thick diagonal $\Delta$ is light blue. Boundary points of $\wideparen{M}_\Delta$
are black, except for those of depth two (shown in dark green).}
\label{fig:lemma}
\end{figure}

On $\wideparen{\Delta}^I_\cC=[0,\infty)\times S^I_\cC$ the 'coordinates' 
$(r^I_\cC,Q^I_\cC)$ are used.
With definitions contained in Appendix \ref{sect:appendix} we get 
the so-called {\em graph blow up}, see \cite[Eq. 8]{AMN}:

\begin{lemma}\quad\\[-6mm] \label{lem:graph:blow:up}
\begin{enumerate}[1.]
\item 
The graph blow up of $M\cong\bR^{(n-1)d}$ by the family 
$\{\Delta^E_\cC\mid\cC\in \cP_\Delta(N)\}$ is 
the topological space
\beq
\wideparen{M}_\Delta := {\rm closure} \big(\hat{I}_\Delta(\widehat{M})\big) \, .
\Leq{wideparen:M:Delta}
It is an $(n-1)d$--dimensional manifold with corners,
see Figure \ref{fig:lemma}). 
\item 
The blow-down map $\beta: \wideparen{M}_\Delta\to M$ is proper, and for binary collisions
\[\beta^{-1}(\Xi_\cC) \ \cong \ \Xi_\cC\times S^I_\cC\qquad (\cC\in  \cP_\Delta(N), |\cC|=n-1).\]
\end{enumerate}

\end{lemma}
\textbf{Proof:}
\begin{enumerate}[1.]
\item 
The family $\{\De_\cC\mid \cC\in \cP_\Delta(N)\}$ is a finite 
semilattice of linear subspaces of $M = \De_{\cC_{\min}}$, that is, closed $p$--submanifolds 
in the sense of Def.\ \ref{def:sub}. As they are linear subspaces, they form a cleanly intersecting
family in the sense of \cite[Def.\ 5.4]{AMN}.
Then it follows from \cite[Theorem 5.12]{AMN} that the   graph blow-up $\wideparen{M}_\Delta$ is a weak
submanifold. Since a weak submanifold of a manifold with corners
is the image of an injective immersion of a manifold with corners, the first statement follows.
\item 
As applied to the present problem, the main statement of \cite[Theorem 5.12]{AMN} is that the 
graph blow-up $\wideparen{M}_\Delta$ is diffeomorphic (in the sense of Definition
\ref{def:smooth}) to the iterated blow-up.
Thus properness of the blow-down map $\beta$ follows from iteration of \cite[Cor.\ 3.7]{AMN}.
Moreover, as the $\Xi_\cC\subseteq\De_\cC$ are the relatively open complements \eqref{Xi:Null},  the 
blow-up of $\Xi_\cC\times \widehat{\Delta}^I_\cC$ at $\Xi_\cC\times\{0\}$  coincides with 
$\Xi_\cC\times\wideparen{\Delta}^I_\cC$. 
\hfill $\Box$\\[2mm]
\end{enumerate}
It follows from \cite[Theorem 5.12]{AMN} that the graph blow up $\wideparen{M}_\Delta$ is diffeomorphic 
to the so-called {\em total boundary blow-up}, see also \cite[Sect.\ 5.13]{Me}.

We identify $\wideparen{M}_\Delta\setminus \partial\wideparen{M}_\Delta$ with $\widehat{M}$.

For some $\alpha\in (0,2)$ we henceforth consider $(-\alpha)$--homogeneous potentials 
$U:\widehat{M}\to \bR$ 
that are of the form \eqref{Ham:Pot} with $Z_{i,j} > 0$ in the two-body potentials 
\[U_{i,j}(q) = Z_{i,j}\, \|q\|^{-\alpha} \quad \big(q\in \bR^d\setminus\{0\}\big) .\] 
In extending the Hamiltonian flow of \eqref{Ham:eq} to the collision manifold, 
the function $U^{-1/\alpha}$ will appear in the differential equations \eqref{qw:p}. 
Extended by zero on $\Delta=M\backslash \widehat{M}$, this is a function $f:M\to \bR$ 
on $M \cong \bR^{(n-1)d}$ that is Lipschitz continuous but not continuously differentiable. 
Instead, when lifting it to the graph blow up $\wideparen{M}_\Delta$, it is in the H\"older space
$C^{(1+\alpha)}(\wideparen{M}_\Delta,\bR)$.%
\footnote{We set $C^{(\alpha)} := C^{k,\alpha'}$ with $k:= \lceil \alpha \rceil-1 \in \bN_0$ 
and $\alpha' := \alpha - k$ for $\alpha \in (0,\infty)$.}
This (together with a similar property of $F$ defined in  \eqref{qw:p}) 
will imply differentiability of the flow. 

Before treating the general case, we will give a simple example.
\begin{example}[boundary defining function]  \label{ex:bdf}
Consider for $\alpha>0$ the function
\beq
f: \bR^2 \to \bR \qmbox{,} f(x) = \left\{\begin{array}{ll}
\big( 1/|x_1|^\alpha + 1/|x_2|^\alpha \big)^{-1/\alpha} &, x\in \widehat{X}\\
0 &, x\in \bR^2\setminus\widehat{X}
\end{array}\right.,
\Leq{Lip:ex}
with $\widehat{X}:=\{x=(x_1,x_2)\in\bR^2\mid x_1\neq0\neq x_2\}$. 
We think of $f$ as the extension of $U^{-1/\alpha}$, with the potential 
of two independent pairs of masses on a line
\[U\in C^\infty(\widehat{X},\bR)\qmbox{,} U(x_1,x_2)= |x_1|^{-\alpha} + |x_2|^{-\alpha} \, .\]
$f$ is 1--homogeneous and
Lipschitz continuous ($f\in C^{0,1}(\bR^2,\bR)$), but $f\not\in C^1(\bR^2,\bR)$,
see Figure \ref{fig:ex}.
\begin{figure}[h]
\centerline{\includegraphics[width=68mm]{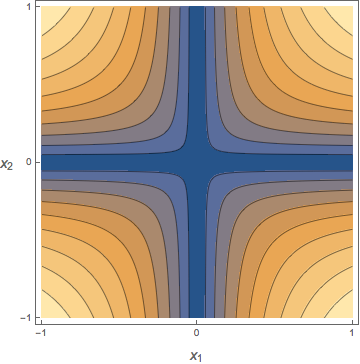}
\includegraphics[width=68mm]{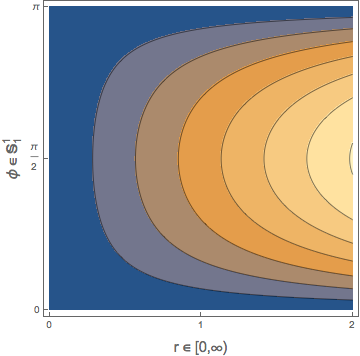}}
\caption{Example \ref{ex:bdf} for $\alpha = 1$: 
The Lipschitz function \eqref{Lip:ex} on $\bR^2$ (left) 
and its lift $\rho\in C^{\infty}(\wideparen{X},\bR)$, for the blow-up $S^1_1\times[0,\infty)$ of one quadrant (right)}
\label{fig:ex}
\end{figure}

\noindent
However, if we continuously extend $U^{-1/\alpha}$
by zero to $\rho:\wideparen{X}\to \bR$ on 
the total boundary blow-up $\wideparen{X}\supseteq\widehat{X}$, then we claim that
$\rho$ is in the H\"older space $C^{(1+\alpha)}(\wideparen{X},\bR)$. 
For the gravitational case $\alpha=1$ we even have
$\rho\in C^\infty(\wideparen{X},\bR)$, and $\rho$ is a boundary defining function 
in the sense of
Melrose.\footnote{
{\bf Definition} ({\sc Melrose} \cite[Lemma 1.6.2]{Me}){\bf :}
A {\bf boundary defining function} on a manifold with corners $X$ is a function $\rho\in C^\infty(X,\bR)$
with $\rho|_{\partial X} = 0$, $\rho|_{X\setminus \partial X} > 0$ and in local coordinates at
$p\in \partial_k X$, $\rho(x) = a(x)x_1\cdot\ldots\cdot x_k$ with $a(p) > 0$ and  $a$ smooth.}
Here $\wideparen{X}$ is diffeomorphic to the disjoint
union of four copies of the manifold with corners 
\[ S^1_1\times[0,\infty)\cong [0,\pi]\times[0,\infty)\qmbox{, with} 
S^m_k:=S^m\cap \bR^{m+1}_k\] 
and half-angle polar coordinates 
$x_1=r \cos(\phi/2)$, $x_2=r\sin(\phi/2)$, $r=\|x\|$. 
Compare with \cite[Lemma 5.10]{AMN} for such pair blow-ups. Then 
\[\rho(\phi,r) = f(x) = \frac{r}{(1/\cos^\alpha(\phi/2)+1/\sin^\alpha(\phi/2))^{1/\alpha}} 
= \frac{r \sin(\phi/2)} {(1+\tan^\alpha(\phi/2))^{1/\alpha}} \, .\]
So $\rho\in C^{(1+\alpha)}(\wideparen{X},\bR)$. For $\alpha =1$, 
$\rho(\phi,r) = \frac{r \sin(\phi)}{\sqrt{8}\sin(\phi/2+\pi/4)}$, proving the claim.
\hfill $\Diamond$
\end{example}
Example \ref{ex:bdf} generalizes as follows:
\begin{lemma}[boundary defining function]\quad\label{lem:rho}\\
The function on the total boundary blow-up, defined 
for $\alpha>0$ by 
\beq
\rho:\wideparen{M}_\Delta\to \bR\qmbox{,}\rho|_{\widehat{M}} := U^{-1/\alpha}\mbox{ and }
\rho|_{\partial\wideparen{M}_\Delta} := 0 \, ,
\Leq{def:rho} 
is in the H\"older space $C^{(1+\alpha)}(\wideparen{M}_\Delta,\bR)$. \\
For the case $\alpha=1$ of celestial mechanics $\rho$ is a boundary defining function.
\end{lemma}
\textbf{Proof:}
Note that by our assumption $Z_{i,j} > 0$, $U > 0$ diverges
to $+\infty$ at $\Delta$. 
It is immediate that $\rho|_{\widehat{M}}>0$ is smooth. 
The behaviour of $\rho$ at the boundary of 
$\wideparen{M}_\Delta$ is given by iterated Taylor expansion with respect to a 
chain for a {\em size order}
of the $\Delta^E_\cC$, see \cite[Section 5.2]{AMN}. For $\cC \in \cP_\Delta(N)$ 
and $q^E_\cC \in \Xi_\cC$ so that  $q^E_C\neq q^E_D$ for $C\neq D\in \cC$ we write 
$q_\cC^I=r Q_\cC^I$ with the notation from \eqref{q:r:Q}. Then we expand
\[ \rho(q) = r \textstyle
\left( \sum\limits_{C\in \cC} \sum\limits_{i<j\in C} \frac{Z_{i,j}}{\|Q_{\cC,i}^I-Q_{\cC,j}^I\|^\alpha}
 + r^\alpha \sum\limits_{C\neq D\in \cC}  \sum\limits_{i\in C,j\in D} 
\frac{Z_{i,j}}{\| q_C^E - q_D^E+ r(Q_{\cC,i}^I-Q_{\cC,j}^I ) \|^\alpha}\right)^{-1/\alpha}\]
with respect to $r$, obtaining a (local) H\"older $C^{(1+\alpha)}$ dependence.\\
For $\alpha=1$ smoothness at the boundary of $\wideparen{M}_\Delta$ 
reduces to Lemma 5.13.3 of \cite{Me}.
\hfill$\Box$
\subsection{Compactifying configuration space} \label{sub:conf:compact}
%
We begin by compactifying the subspaces $ \Delta^E_\cC\subseteq M$ 
in $\overline M= M \sqcup \bS$,
resulting in closed disks $\overline{\Delta}^E_\cC$ whose boundary 
$\bS\cap \overline{ \Delta}^E_\cC$ is a subsphere of $\bS$.\\
With $\cP_\bS (N) = \cP (N)\backslash \{\cC_{\max}\}$ from \eqref{cP:Delta}, 
modifying \eqref{hat:I:Delta}, we use the map
\[\hat{I}:\widehat{M} \longrightarrow  
\prod_{\cC\in \cP_\Delta (N)} \!\!\! \big[\overline{M}:\overline{ \Delta}^E_\cC\big] \
\times \prod_{\cC\in \cP_\bS (N)} \!\!\!  
\big[\overline{M}: (\bS\cap \overline{ \Delta}^E_\cC)\big] 
\ \mbox{,}\quad
q\longmapsto (q)_{\cC\in \cP_\Delta(N) \sqcup \cP_\bS(N)},\]
which smoothly imbeds $\widehat{M}$ as a submanifold. We omitted $\cC_{\max}$
in the second product,
since $\bS\cap \overline{ \Delta}^E_{\cC_{\max}} =\bS\cap \{0\} =\emptyset$.\\
Note that the $\overline{ \Delta}^E_\cC$ intersect the boundary $\bS$ of $\overline{M}$
neatly (see, {\em e.g.} {\sc Hirsch} \cite[Section 1.4]{Hi}), and are p--submanifolds
of $\overline{M}$ in the sense of  \cite[Definition 1.7.4]{Me}. 
So like in Lemma \ref{lem:graph:blow:up} the 
graph blow up 
\beq
\wideparen{M} := {\rm closure} \big(\hat{I}(\widehat{M})\big) 
\Leq{wideparen:M}
of the compact ball $\overline{M}$ has the structure of an 
$(n-1)d$--dimensional manifold with corners. The new feature is that, 
unlike $\wideparen{M}_\Delta \subseteq \wideparen{M}$ defined in 
\eqref{wideparen:M:Delta} and
$\widehat{M}\sqcup \widehat{\bS} \subseteq \wideparen{M}$ from \eqref{wideparen:M:S}, 
$\wideparen{M}$ is compact.  
\begin{figure}[h]
\centerline{\includegraphics[width=80mm]{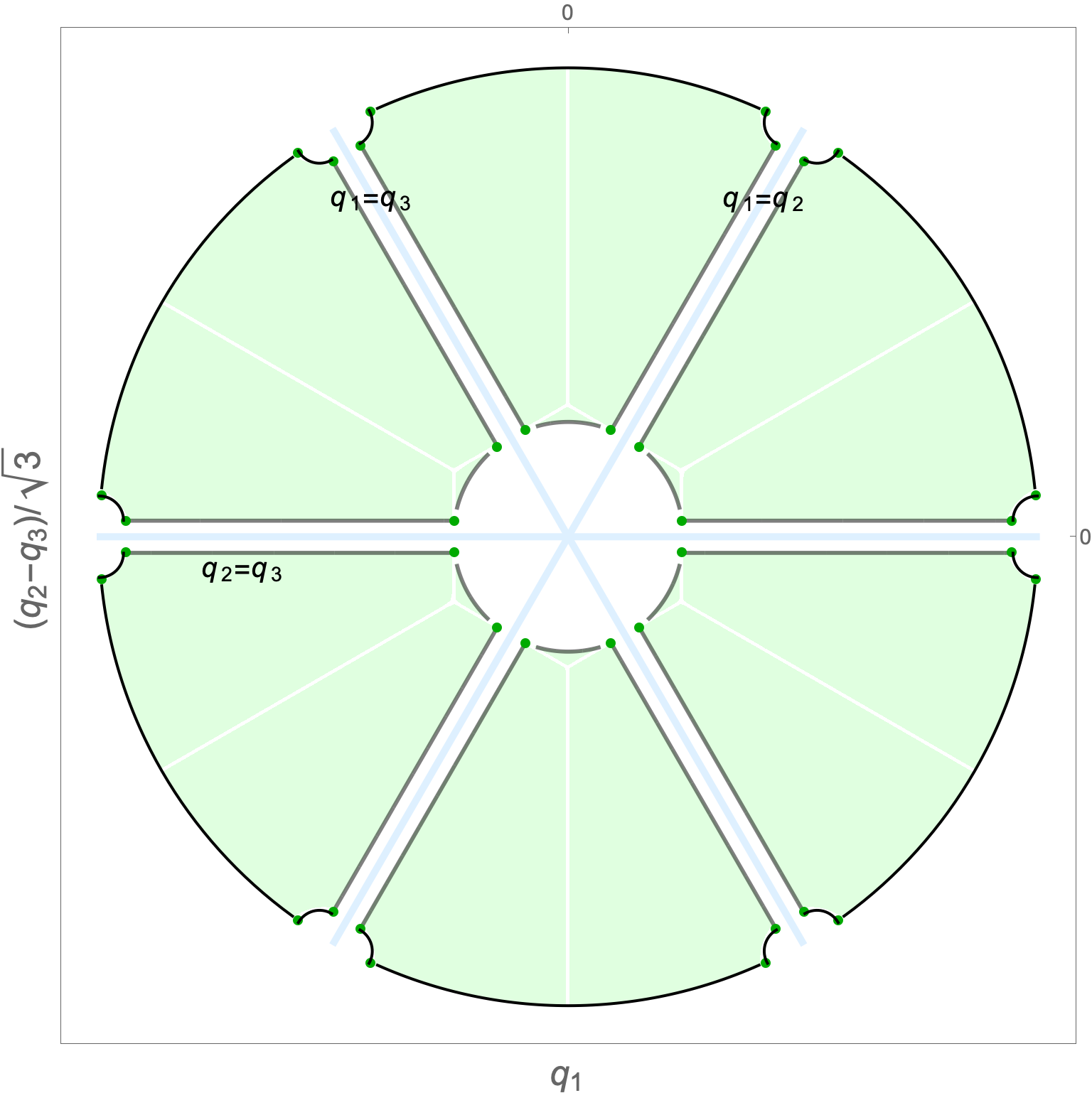}}
\caption{Schematic view of the blown-up configuration space $\wideparen{M}$
(light green)
arising by blowing up configuration space 
$M \cong \bR^2$ for $n=3$ particles in $d=1$ dimension. 
The thick diagonal $\Delta$ is light blue. Boundary points of $\wideparen{M}$
are black, except for those of depth two (shown in dark green).}
\label{fig:blownupConfigSpace}
\end{figure}
\section{Blowing up the energy surfaces}\label{sec:blow:energy}
%
We will use Section \ref{sub:blow:conf} to blow up and compactify  the energy surfaces 
$\hSE = H^{-1}(E)$. 
The  resulting  compact manifolds with corners,  denoted   $\wideparen{\Sigma}_E$,  are easy to define
and understand via their projection to  $\wideparen{M}$.   See equation \eqref{sphere:bundle}. %

Analysing the dynamics induced on  the  
flow-invariant boundary components takes some work.
After the introductory example of two bodies in Subsection \ref{sub:two:bodies},
we separately consider the  dynamics on the pieces of the energy surfaces 
over the boundary blow-ups $\wideparen{M}_{\widehat \bS}$, $\wideparen{M}_\Delta$ 
and finally over the entire  compact 
total blow-up $\wideparen{M}$ of configuration space 
$\widehat{M}: = M \setminus \Delta$.

\subsection{Rescaling velocities} 
\label{subsec: rescale vel}
We prefer to work with velocities 
$$v = \cM^{-1}\, p$$
rather than momenta.  Then the energy is given by
$$E =  \textstyle{\frac{1}{2}} \| v \|_{\cM}^2 - U(q)$$
Make the position dependent rescaling of velocities
\beq
w := G_E(q)v  \quad \mbox{ with } \textstyle
G_E :=  \big( 2(E+ U)\big)^{-\eh} \, ,
\Leq{w:tau}
so that the energy becomes
$$E = \big(E + U(q)\big) \|w \|^2 - U(q) \, .$$
Solving, we see that the energy in the $(q,w)$--variables  is $E$ iff $\|w \| = 1$.

The substitution (\ref{w:tau}) requires that $E + U \ge 0$
(we don't want imaginary $w$'s!) which means that $q$ 
must lie in the Hill region
 \[ \widehat{M}_E := \big\{q\in \widehat{M}\mid U(q)\ge -E\big\} \, . \] 
 If $E \ge 0$ this is no restriction since our $U$ is positive everywhere
 and so $\widehat{M}_E = \widehat M$. 
 In this case the energy level set $\hSE$ is equal to $\widehat M \times \bS$
 where $\bS= {\mathbb S}^{(n-1)d-1}$ is the unit sphere in the $w$--variables.   
 We can now simply take closures by letting $q \to \partial \widehat M$
 and realizing our spheres are staying constant.
 So for $E \ge 0$ we arrive at  
 \[\wideparen{\Sigma}_E  =  \wideparen{M} \times \bS \, .\]  
When $E < 0$  we must pay  attention to the behaviour of velocities 
as we approach  the Hill boundary $\{ U = -E \}$.
In the interior of the Hill boundary the fibers of the original velocity projection 
$(q,v) \mapsto q$
restricted to $\hSE$ are spheres whose radius shrinks so that they  
degenerate to points when we reach the Hill boundary.
In our new $w$--variables the spheres are all the same size, 
so as we approach the boundary they remain the same.
In going to the $w$--variables we thus   replace the original smooth $\hSE$ by 
$\widehat{M}_E \times \bS$.
The extraneous directions along the boundary will not cause problems 
with the dynamics, but will 
require a bit of analysis in order to make sense of brake orbits -- 
orbits hitting the Hill boundary.  This is done in
Subsection \ref{subsec: Hill boundary analysis} .      
Letting $q$ tend to boundary points of $\wideparen{M}$ while
remaining in $\widehat{M}_E$ we see that  again
the $w$--spheres do not change size:  we just have $w \in S$, the fixed unit sphere.  
Thus regardless of the energy $E$ we get  
\beq
 \wideparen{\Sigma}_E \cong  \wideparen{M}_E \times \bS
\Leq{sphere:bundle}
where  $\wideparen{M}_E$ denotes the closure of the smooth manifold with boundary
$\widehat{M}_E$ within~$\wideparen{M}$. 

\begin{remarks}\qquad\\[-6mm] 
\begin{enumerate}[1.]
\item 
{\bf Rescaling of velocities:}
Our rescaling of velocities is half of our globalization of the construction first devised by 
{\sc McGehee} in \cite{McG1}, and developed further
by {\sc Devaney} \cite{Dev} and many others.
See in particular {\sc Lacomba} and {\sc Ibort} \cite{LI}. 
The other  half comes with the next section when  we rescale time.
\item 
{\bf Critical points:}
For our choice of $U$ the function $U$ has no critical values and hence
the Hill boundary is smooth regardless of the choice of energy $E$.
\item 
{\bf Tangent space:} 
Note here that the tangent space $TX$ of a manifold with corners $X$
like $\wideparen{M}_E$ is well-defined, since 
$X$ is contained in a smooth manifold $\widetilde{X}$ of the same dimension. 
$TX$ can then be defined invariantly by the restriction 
of the bundle $T\widetilde{X}\to \widetilde{X}$
to $X$, see \cite[Remark 2.7]{AMN}.
\hfill $\Diamond$
\end{enumerate}
\end{remarks}
\subsection{Rescaling time} 
\label{subsec: rescale time}
To extend the differential equations to $\wideparen{\Sigma}_E$, we introduce a new 
time parameter $\tau$ along orbits with 
\beq
{\textstyle \frac{dt}{d\tau} = \, \widetilde{G}_E(q)\, G_E(q) \quad
\mbox{for } \widetilde{G}_E :=\frac{E+U}{1+U} U^{-1/\alpha}\, }
\Leq{new:time}
and denote $\frac{d}{d\tau}$ by '.
If we set
\beq \textstyle
F(q) := \frac{{\cal M}^{-1}\nabla U(q)} {2(1+U(q))\,U(q)^{\frac{1}{\alpha}}} \,,
\Leq{force:term}
the Hamiltonian equations \eqref{Ham:eq} on $\widehat{P}$ acquire the form
\beq
q' =  \widetilde{G}_E(q) \, w \qmbox{,} w' =   F(q) - \langle F(q),w\rangle\, w  
\Leq{qw:p}
leaving the energy surfaces $\hSE$ invariant.  
As $q\to0$, the force term $F(q)$ is asymptotically homogeneous of degree zero. 
It is bounded, since the terms 
\beq 
\frac{\nabla_{\!q} U_{i,j}(q)} {U_{i,j}(q)^{1+\frac{1}{\alpha}}} 
= - { \alpha Z_{i,j}^{-\frac{1}{\alpha}} \, \frac{q}{\|q\|} }
\qquad \big(q\in \bR^d\setminus\{0\}\big)
\Leq{two:body:hom}
are bounded and since $U_{i,j} > 0$.
\begin{remarks}[scalings] \quad\\[-6mm]
\begin{enumerate}[1.]
\item 
{\bf Change of speed:} To normalize it to one is just a convenient choice.
\item 
{\bf Time change:}
The factor $U^{-1/\alpha}$ in the choice \eqref{new:time} 
of the time change is motivated by the desire that the velocity 
$q'$ should vanish asymptotically linearly in the distance from the boundary.
Looking at the first equation in \eqref{qw:p}, this is more or less obvious for boundary points
in $\Delta$. 
As we shall see, this is also true for  boundary points in the sphere $\bS$ at infinity.
\item 
{\bf Alternatives:}
One drawback of our choice is that $\|w\|=1$ also 
at the boundary $\partial \widehat{M}_E$ of Hill's region.\\
Therefore in our companion paper \cite{KM} we use scalings for speed and time change that 
depend on total energy $E$.
\hfill $\Diamond$
\end{enumerate}
\end{remarks}
Before looking carefully   at the different aspects of the compactification,
we present the simplest example.

\subsection{The example of two bodies}\label{sub:two:bodies}
%
The simplest case is the reduced system with Hamiltonian $\|p\|^2/2-U(q)$ and
$U(q) = \|q\|^{-\alpha}$ where $q \in \bR^d$   stands for $q_1 -q_2$. 
 Hamiltonian equations are
\beq
\dot{q}=p\qmbox{,}\dot{p} = \nabla U(q)
\Leq{2body:ham:ode}
with $\nabla U(q) = - \alpha\, Q/\|q\|^{\alpha+1}$ for $Q := q/\|q\|$, 
so that our force term \eqref{force:term} equals
\[\textstyle F(q) = - \frac{\alpha}{2}\frac{ Q}{1+\|q\|^\alpha} \,. \]
The scaling functions \eqref{w:tau} for velocity and \eqref{new:time} for time
have the form 
\[\textstyle
G(q) = \big( 2(E + \|q\|^{-\alpha}) \big)^{-1/2} \qmbox{,}
\widetilde{G}_E(q) =\frac{E+ \|q\|^{-\alpha}}{1+ \|q\|^{-\alpha}} \|q\| \,.\]
For total energy $E\in \bR$ the differential equation \eqref{qw:p}
extends to the boundary of the energy surface 
$\wideparen{\Sigma}_E = \wideparen{M}_E\times \bS^{d-1}$ over compactified
configuration space
\[\wideparen{M}_E:= \left\{
\begin{array}{ccl}
[0,\infty] \times \bS^{d-1} &,& E\ge0\\
\big[0,|E|^{-1/\alpha} \big] \times \bS^{d-1} &,& E<0
\end{array}
\right.\ .\]
The boundary of $\wideparen{\Sigma}_E$   has two components,   one with $r = 0$
and the other with either $r = \infty$ or $r = |E|^{-1/\alpha}$ where $r = \|q \|$.

\begin{figure}[h]
\centerline{\includegraphics[width=44mm]{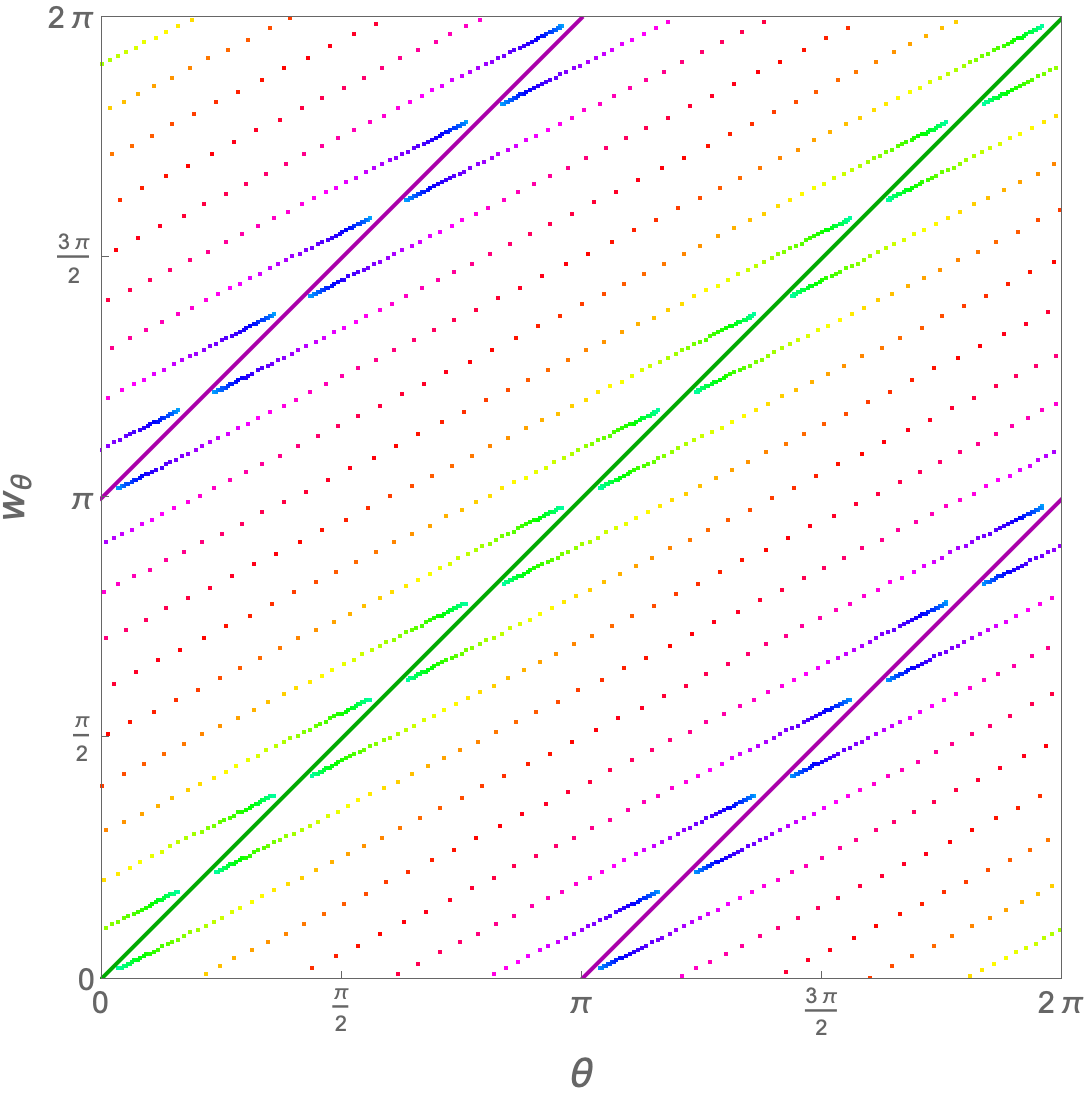}
\hfill
\includegraphics[width=44mm]{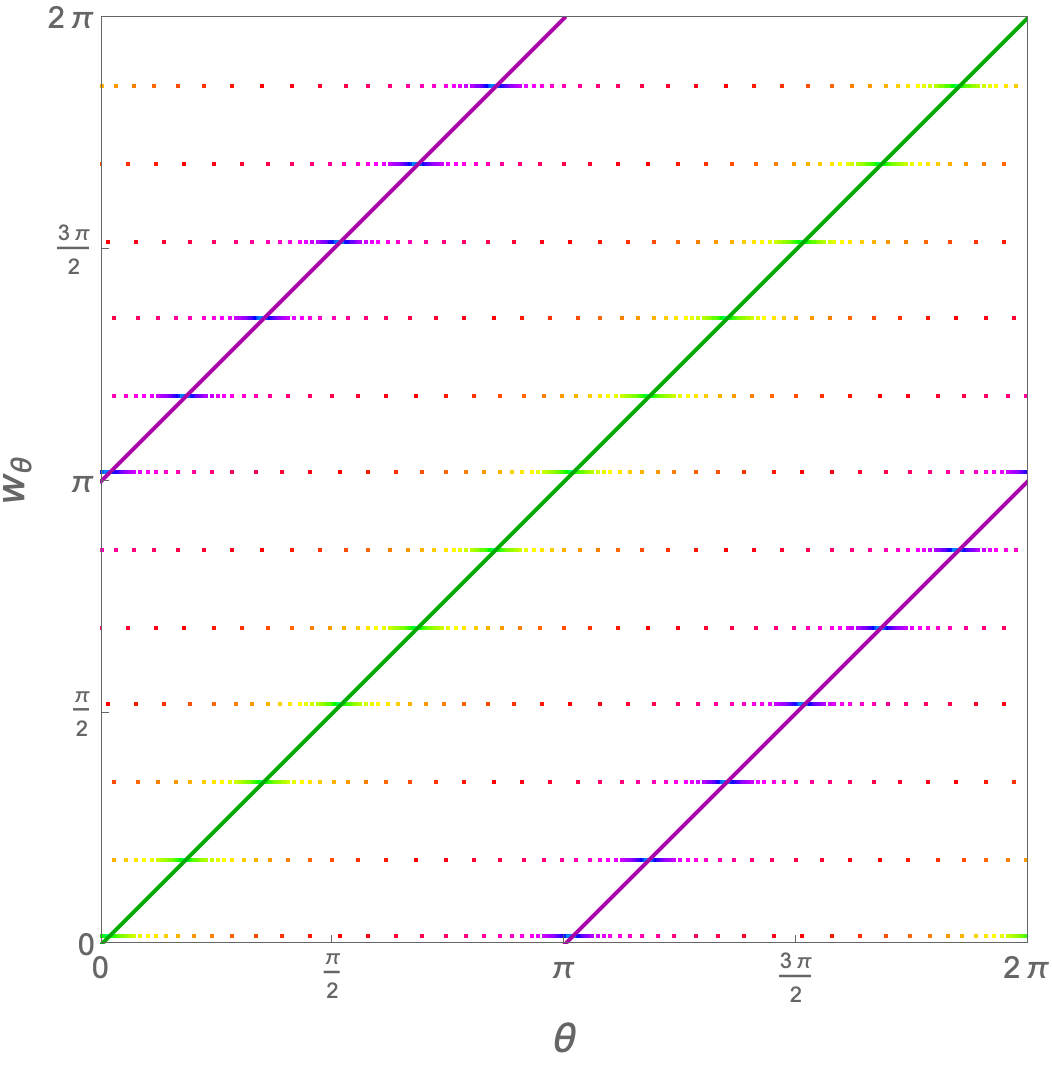}
\hfill
\includegraphics[width=44mm]{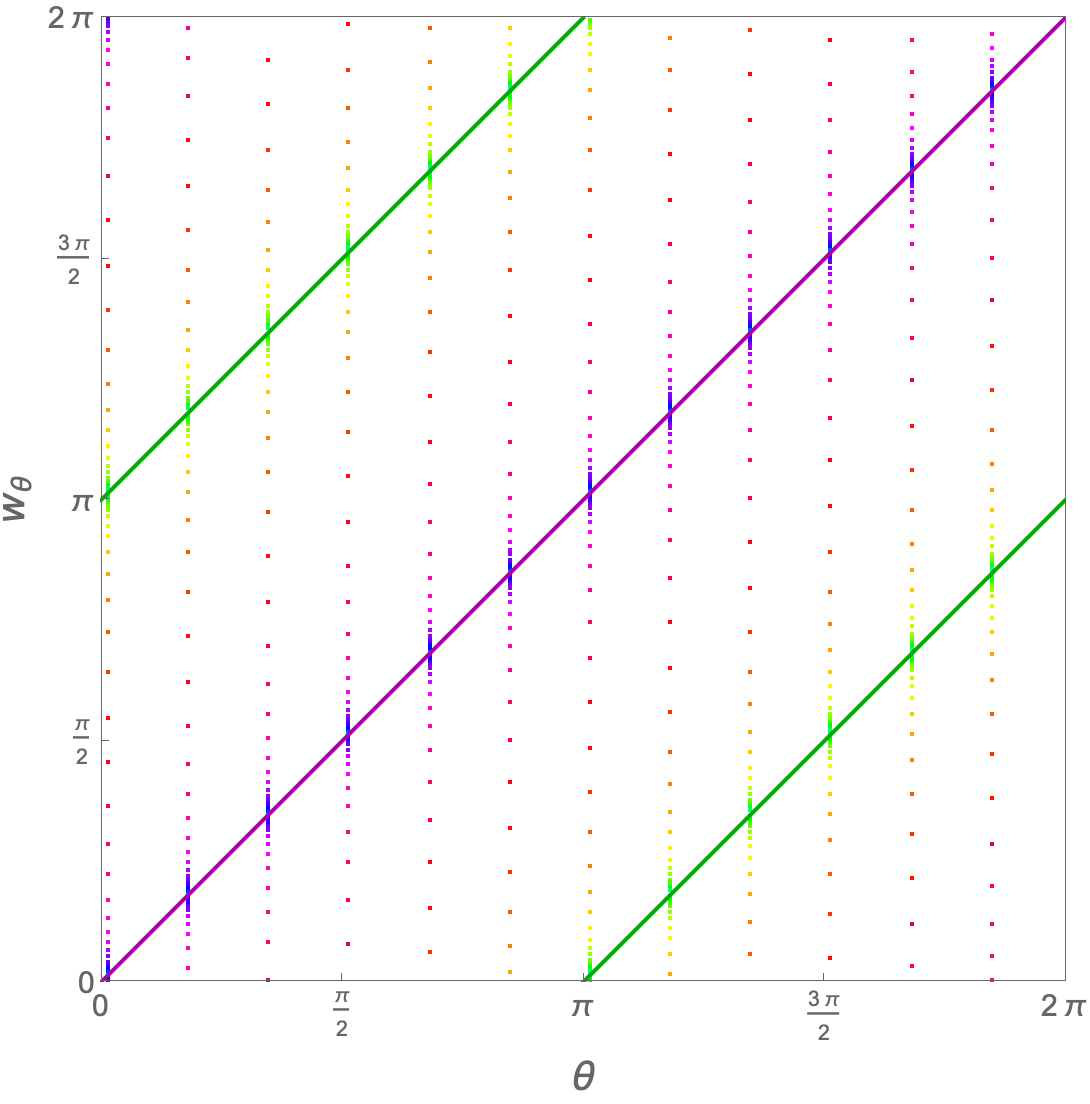}}
\caption{Flows on invariant boundary tori, for the gravitational case $\alpha=1$. 
Left: At collision. Middle: At infinity, for energy $E=1$. 
Right: At the boundary of Hill's region, for energy $E=-1$. The unstable rest points
are coloured in magenta, the stable ones in green. The flow lines are coloured by time 
$\tau$.} 
\label{fig:boundaries}
\end{figure}

Using $Q=q/r$ one gets 
\beq \textstyle
r'= r \frac{1+Er^\alpha}{1+r^\alpha} \langle Q,w\rangle 
\ \mbox{,}\
Q' = \frac{1+Er^\alpha}{1+r^\alpha} \big(w- \langle Q,w\rangle Q \big)
\ \mbox{,}\
w' = \frac{-1}{1+r^\alpha} \frac{\alpha}{2}  \big(Q- \langle w,Q\rangle w \big)\,.
\Leq{rQw}
As $\mathrm{span}(Q,w)$ is invariant, it suffices to consider dimension $d=2$.
So we use polar coordinates for the $q$ and $w$ variables, 
assumed to be complex-valued:
\[q = r \exp(\imath\,\theta) \qmbox{,} w = \exp(\imath\, w_\theta) \, .\]
Setting
$$\psi = w_\theta-\theta$$
the   differential equation takes the form
\beq \textstyle
r'= r \frac{1+Er^\alpha}{1+r^\alpha} \cos(\psi) \qmbox{,} \theta' = \frac{1+Er^\alpha}{1+r^\alpha} \sin(\psi)
\qmbox{,}
w_\theta' = \frac{\alpha}{2} \frac{1}{1+r^\alpha} \sin(\psi) \, .
\Leq{d:e:2body}
The $r$ equation shows that the two boundary components  are invariant.  They are 
tori coordinatized by $(\theta, w_{\theta})$.  On these tori the circles $\{\theta = w_\theta\}$ and 
$\{\theta = w_\theta+\pi\}$ consist of rest points (see Figure \ref{fig:boundaries}).

\begin{enumerate}[$\bullet$]
\item 
{\bf Collision:}\\
At the flow invariant boundary component of $\wideparen{\Sigma}_E$ with 
$q=0$ the o.d.e.\  \eqref{d:e:2body} equals
\[
 \textstyle
\theta' = \sin(w_\theta-\theta)
\qmbox{,}
w_\theta' = \frac{\alpha}{2} \sin(w_\theta-\theta) \, .
\]
So the differential equation $\psi' = -(1-\frac{\alpha}{2})\sin(\psi)$   
or $\psi:=w_\theta-\theta $ is solved for initial conditions with 
$w_\theta(0) = \theta(0)\pm \pi/2$
by $\psi(\tau)= \pm 2 \cot ^{-1}\left(e^{(1-\frac{\alpha}{2})\tau }\right)$,
whereas  $w_\theta - \frac{\alpha}{2}\theta  = (1-\frac{\alpha}{2})\theta(0)\pm \pi/2$ 
is a constant of the motion. We obtain 
\[ \textstyle
\theta(\tau) 
= \theta(0) \pm
\frac{\frac{\pi}{2} - 2 \cot ^{-1}(e^{(1-\frac{\alpha}{2})\tau})}
{1-\frac{\alpha }{2}}\ \mbox{ , }\
 w_ \theta(\tau) =  \theta(0) \pm
 \frac{ \pi/2 - \alpha \cot ^{-1}\left(e^{(1-\frac{\alpha}{2})\tau }\right)}{1-\frac{\alpha}{2}}
 \, .\]
The total angle covered is  
\[  \lim_{\tau\nearrow \infty} \textstyle \big(\theta(\tau)-\theta(-\tau) \big) 
= \pm\frac{\pi}{1-\frac{\alpha}{2}} \, .\]
This agrees, as it should,  with the range of `Rutherford' type 
scattering for homogeneous central forces as calculated in \cite[Section 4]{KK}.
The limit points 
\begin{align*} 
\lim_{\tau\to +\infty}  \big(\theta(\tau), w(\tau)\big) 
&= \textstyle{\big(\theta(0)\pm \frac{\pi \, }{2-\alpha}\, , \,
\theta(0)\pm \frac{\pi }{2-\alpha} \big)} \, ,\\
\lim_{\tau\to - \infty}  \big(\theta(\tau), w(\tau)\big) 
&= \textstyle{\big(\theta(0)\mp \frac{\pi \, }{2-\alpha} \, , \,
\theta(0)\mp (\frac{\pi }{2-\alpha} - \pi) \big)} 
\end{align*} 
are stable respectively unstable rest points. So the two unstable orbits of a rest point
converge to the {\em same} stable rest point exactly in the cases
\[\alpha =2(1-1/m)\qquad \big( m\in\bN\,\backslash\{1\} \big).\] 
Only then we can uniquely regularize the original o.d.e.\
\eqref{2body:ham:ode} at collision.
The cases $m$ odd then
correspond to motion in the forward direction,
whereas $m$ even (including the gravitational case $\alpha=1$) 
corresponds to backscattering.
\item 
{\bf Spatial infinity:}\\
Similarly, for $E\ge0$ the boundary component at $r=\infty$ is invariant under the flow,
as one sees by considering the o.d.e. 
$\rho' = -\rho \frac{E+\rho^\alpha}{1+\rho^\alpha} \langle Q,w\rangle Q$ for $\rho := 1/r$.
Thus at spatial infinity \eqref{d:e:2body} takes the form
\[
\theta' = E \sin(w_\theta-\theta)
\qmbox{,}
w_\theta' = 0 \, ,
\]
with the solutions for $w_\theta(0) = \theta(0)\pm \pi/2$
\[  \textstyle
\theta(\tau) = \theta(0)  \pm \big(\frac{\pi}{2} - 2\tan^{-1}\left( \exp(-E \tau) \right)
 \big)
 \qmbox{,}
w_ \theta(\tau) = w_ \theta(0) \, .\]
They converge to the rest points 
$\big(\theta(0)\pm \pi/2, w_ \theta(0)\big) $ as $\tau\nearrow\infty$ and\\ 
$\big(\theta(0)\mp \pi/2, w_ \theta(0)\big) $ as  $\tau\searrow-\infty$ .
So the stable manifold of a rest point $(\theta,\theta)$ equals the unstable manifold of
$(\theta+\pi,\theta)$, see Figure \ref{fig:boundaries}, middle.\\
For energy $E=0$ the whole boundary $\bS^1\times \bS^1$
at infinity consists of rest points.
All energy $E$ solutions of \eqref{rQw} (except those colliding) thus have the property
that they go to a stable fixed point as $\tau\nearrow \infty$.
\item 
{\bf Boundary of Hill's region:} \\
For energy $E<0$ the boundary
$\partial \widehat{M}_E= \{ q\in \bR^2\mid \|q\| = |E]^{-1/\alpha} \}$
of Hill's region is invariant under the flow, too.
Then \eqref{d:e:2body} becomes
\[\textstyle
\theta' =0  \qmbox{,} 
w_\theta' = \frac{\alpha}{2} \frac{1}{1+1/|E|} \sin(w_\theta-\theta) \,,\]
with the solutions for $w_\theta(0) = \theta(0)\pm \pi/2$ 
and $c:= \frac{\alpha}{2} \frac{1}{1+1/|E|} $
\[  \textstyle
\theta(\tau) = \theta(0) 
 \qmbox{,} w_ \theta(\tau) = 
w_ \theta(0)  \pm \big(\frac{\pi}{2} - 2\tan^{-1}\left( \exp(c \, \tau) \right)\big)\, .\]
Here the unstable manifold of a rest point $(\theta,\theta) \in \bS^1\times \bS^1$ 
equals the stable manifold of
$(\theta,\theta+\pi)$, see Figure \ref{fig:boundaries}, right.
As $(\theta,\theta)$ and  $(\theta,\theta+\pi)$ are uniquely connected in the positive time
direction,
one can uniquely connect the incoming and outgoing solutions
(brake orbits) for the original o.d.e.\ \eqref{2body:ham:ode}.
(See also the end of Subsection \ref{dynamic philosophy} below.) 
\end{enumerate}
\subsection{Blowing up the energy surface at Hill's boundary}
\label{subsec: Hill boundary analysis}
%
The Hill boundary 
$\partial \widehat{M}_E\subseteq  \widehat{M}$ is a submanifold of codimension one.
(All values $E$ of our potential are regular values.) 
So (like in Subsection \ref{sub:conf:infty}) 
blowing up Hill's  region $\widehat{M}_E$ at $\partial \widehat{M}_E$ just reproduces 
$\widehat{M}_E$. The restriction of the trivial sphere bundle \eqref{sphere:bundle}
to $\partial \widehat{M}_E\subseteq \partial\wideparen{M}_E$ has the form
\[\partial \widehat{M}_E \times \bS^{(n-1)d-1} \to \partial \widehat{M}_E\, .\]
\begin{lemma}[Hill's boundary]\quad\\  
For all energy values $E<0$ the flow on $\wideparen{\Sigma}_E$ leaves the fibers 
\[ (\pi_E)^{-1}(q_0) \subseteq \wideparen{\Sigma}_E  
\qquad \big( q_0 \in \partial \widehat{M}_E \big) \]
invariant. The limit point $(q_0,w)$ of the incoming brake orbit is
uniquely connected by its unstable manifold 
to the limit point $(q_0,-w)$ of the outgoing brake orbit.
\end{lemma}
\textbf{Proof:}
The sphere $(\pi_E)^{-1}(q_0)$ is invariant under the
flow generated by the differential equation \eqref{qw:p}, since  $\widetilde{G}_E(q)=0$
for $q_0 \in\partial \widehat{M}_E$, so that $q'=0$.  The vector field
$-\nabla U(q)\neq0$ is an outward pointing normal to $T_q\partial \widehat{M}_E$. 
As $w' =   F(q) - \langle F(q),w\rangle\, w$, the only rest points $(q_0,w)$ 
of the flow on the sphere are the ones with $w$
parallel or antiparallel to the force $F(q_0) = \nabla U(q_0)$.   But the brake orbits
have these incoming and outgoing directions as is seen by a Taylor expansion
of Newton's equations $\ddot q = \nabla U(q)$ with $\dot q(0) = 0, q(0) = q_0$.
 \hfill $\Box$
\begin{remark} [significance of the flow on $\wideparen{\Sigma}_E$]%
\quad\label{dynamic philosophy}\\
We take the general  philosophy that one can concatenate two solution curves to the
reparameterized flow on $\wideparen{\Sigma}_E$ if 
\begin{enumerate}[1.]
\item 
the $\omega$--limit set of first is the $\alpha$--limit set of the second,
with their common limit point $R\in\wideparen{\Sigma}_E$ being a rest point of 
the extended flow.
Note that all the rest points are in $\partial \wideparen{\Sigma}_E$.
\item 
Additionally, we demand that $R$ does not belong to the $\omega$--
respectively $\alpha$--limit sets of other solution curves (disregarding the constant solution
$R$).
\end{enumerate}
In this way we regain
the standard brake orbit solutions for Newton's equations as follows.  
Let $q_0 \in  \partial\widehat{M}_E$ be a brake point
and consider the energy $E$ Newtonian   solution $q(t)$ for which $q(0) = q_0$ 
and consequently $\dot q (0) = 0$.  This solution satisfies $q(-t) = q(t)$.
In the reparameterized time and with rescaled velocities $w$
this single brake curve blows up into the concatenation 
$\gamma_- * \gamma_0 * \gamma_+$ of three curves. 
The $\pi_E$ projections of the two  curves $\gamma_{\pm}$ 
lie in the interior of the Hill region where they  are reparameterizations
of $q(t)$,  parameterized so that they approach $q_0$ as $t \to \pm \infty$.  
The middle curve $\gamma_0$
travels along the invariant sphere   $(\pi_E)^{-1}(q_0)$, connecting the incoming
normalized velocity   
$w = -\nabla U(q_0)/\| \nabla U (q_0)\|$ (the outward pointing  normal to the Hill boundary)  
to the outgoing normalized velocity
$w = +\nabla U(q_0)/\| \nabla U (q_0)\|$, taking infinite  $\tau$-time to do so.
\hfill $\Diamond$
\end{remark}
\subsection{Blowing up the energy surface at collisions}\label{sub:energy:coll}
%
We next consider the boundary component of $\wideparen{\Sigma}_E$
defined in \eqref{sphere:bundle} projecting to $\Delta$ and thus set
\[\wideparen{\Sigma}_{E,\Delta}:= \big\{x\in \wideparen{\Sigma}_E \mid \pi_E(x)\in M \big\}\, .\]
That is, we defer the analysis of the boundary component at spatial infinity.\\
The blown up energy surface $\wideparen{\Sigma}_{E,\Delta}$ has two types of
boundary components, the one projecting to the collision set $\Delta$ and, for $E<0$,
the ones projecting to the boundary of Hill's region. 
Unlike in the last subsection, over $\Delta$ the vector field experiences a loss of smoothness:
\begin{lemma}\quad\\  
The smooth vector field $\widehat{X}_E: \widehat{\Sigma}_E\to T\widehat{\Sigma}_E$
defined by the right hand sides of \eqref{qw:p} continuously extends to a locally
$C^{(1+\alpha)}$ H\"older continuous vector field 
\[\wideparen{X}_{E,\Delta}: 
\wideparen{\Sigma}_{E,\Delta} \to T\,\wideparen{\Sigma}_{E,\Delta}\,. \]
Its flow leaves the boundary component of $\wideparen{\Sigma}_{E,\Delta}$ over $\Delta$
invariant.
\end{lemma}
\textbf{Proof:}
The factor $\frac{E+U}{1+U}$ of $\widetilde{G}_E$ 
has the constant limit 1 over $\Delta$. The factor
 $U^{-\frac{1}{\alpha}} $ of $\widetilde{G}_E$ has been shown
in Lemma \ref{lem:rho} to extend to
a $C^{(1+\alpha)}$ H\"older continuous function on the blow-up $\wideparen{M}_\Delta$.
Since it goes to zero at $\Delta$, $\widetilde{G}_E$ extends to a $C^{(1+\alpha)}$ function
on the blown up configuration space $\wideparen{M}_\Delta$ 
(see \eqref{wideparen:M:Delta}), vanishing over $\Delta$.
So by the first differential equation in \eqref{qw:p} the flow leaves 
the boundary component of $\wideparen{\Sigma}_{E,\Delta}$ over $\Delta$ invariant.

The argument for the force terms $F$ 
is similar, using that the radial blow-up of \eqref{two:body:hom}
is smooth.
\hfill $\Box$
\begin{corollary}[smoothness of the flow]\quad\\
The initial value problem 
$x'= \wideparen{X}(x)$, $x(0)=x_0\in \wideparen{\Sigma}_{E,\Delta}$,
derived from \eqref{qw:p} has unique local solutions 
in $C^1(D,\wideparen{\Sigma}_{E,\Delta})$, with 
open domain $D\subseteq \bR_\tau \times \wideparen{\Sigma}_{E,\Delta}$ containing 
$\{0\}\times  \wideparen{\Sigma}_{E,\Delta}$.
\end{corollary}
\textbf{Proof:}\\
The vector field $\wideparen{X}_{E,\Delta}$ is in the H\"older space 
$C^{(1+\alpha)} \big( \wideparen{\Sigma}_{E,\Delta} , T\,\wideparen{\Sigma}_{E,\Delta} \big)$
with $\alpha>0$ and thus fulfills the criterion of the theorem of  Picard and Lindel\"of.
\hfill $\Box$
\begin{remark}[smoothness of the flow]\quad\\ 
We believe that the flow is $C^{(1+\alpha)}$ H\"older continuous, since the vector field
is $C^{(1+\alpha)}$ H\"older continuous. However, we couldn't find such a result in the 
literature, and we didn't try to prove it.  
\hfill $\Diamond$
\end{remark}

We now consider the flow at collisions more precisely. Every $\cC\in \cP_\Delta(N)$
defines a different type of how the particles meet, as
the $\Xi_\cC$ defined in \eqref{Xi:Null} lead to the stratification 
\beq 
\Delta = \bigsqcup_{\cC\in \cP_\Delta(N)}\ \Xi_\cC
\Leq{Delta:decomposition}
of the thick diagonal $\Delta$.
For $q_0\in \Xi_\cC\subseteq \Delta^E_\cC$, we use the 
$\cM$--orthogonal decomposition 
$M = \Delta^E_\cC\oplus \Delta^I_\cC$ of configuration space, with dimensions 
given in \eqref{dim:MC}. As $\Xi_\cC\subseteq \Delta^E_\cC$ is relatively open, we can 
use the local Cartesian coordinates $(q^E_\cC,q^I_\cC)$ from \eqref{cl:bar} in a neighborhood
$U_\cC:= U^E_\cC\times U^I_\cC\subseteq \Xi_\cC\times \Delta^I_\cC$ of $q_0$. 
Depending on the choice of  $U^E_\cC$, we can choose $U^I_\cC$ so that 
\[U_\cC\cap \Delta= U^E_\cC\qmbox{and} 
U_\cC\setminus U^E_\cC \subseteq \widehat{M}_E^{\mathrm{int}} \, .\]
With the definition \eqref{J:decomp} of $J^I_\cC$ and regarding it as a function on 
$\Delta^I_\cC$, we write the coordinate $q^I_\cC$ 
in the form
\beq
q^I_\cC = r \,Q^I_\cC \quad\mbox{with } J^I_\cC(Q^I_\cC)=1 \mbox{ and } 
r:=\big(J^I_\cC(q^I_\cC) \big)^{1/2}.
\Leq{def:r:Q}
Assuming that $q^I_\cC\neq 0$, this polar decomposition is unique and
\[\Delta^I_\cC\setminus\{0\} \cong (0,\infty)\times S^I_\cC \quad \mbox{ with the sphere }
S^I_\cC:=\big(J^I_\cC\big)^{-1}(1) \, .\]
For any $q_0\in \Xi_\cC$ and any direction $Q^I_\cC\in S^I_\cC$ for which 
$q_0+rQ^I_\cC\not\in \Delta$ if $r>0$ is small enough, 
we thus attach a point to blow up the boundary. 
The last condition is only violated for the subset 
$S^I_\cC\cap \Delta$, which by \eqref{Delta:decomposition}
is a finite union of submanifolds of codimension at least $d$ in $S^I_\cC$.

Although the potential $U$ is in general not $(-\alpha)$--homogeneous in $r$ for $q=q_0+rQ^I_\cC$, 
it has this property asymptotically as $r\searrow 0$: 
For $Q^I_\cC \in S^I_\cC\backslash \Delta$
\[W_\cC(Q^I_\cC) := \lim_{r\searrow 0} r^\alpha\, U(q_0+rQ^I_\cC)
= \eh\sum_{C\in \cC} \sum_{i\neq j\in C} \! U_{i,j} \big( (Q^I_\cC)_i - (Q^I_\cC)_j \big) 
= U^I_\cC (Q^I_\cC)\]
has values in $(0,\infty)$,
so $W_\cC= U^I_\cC\big|_{S^I_\cC}$ and
\beq
\nabla W_\cC(Q^I_\cC) = \nabla U^I_\cC(Q^I_\cC)  - 
\langle U^I_\cC(Q^I_\cC), Q^I_\cC\rangle Q^I_\cC \, .
\Leq{W:C:p}
This limit does not depend on $q_0\in \Xi_\cC\subseteq \Delta^E_\cC$ 
and defines a smooth function 
\[W_\cC:S^I_\cC\setminus \Delta\to (0,\infty) \,.\]
To construct our collision manifold, we attached to the point $(q_0,Q^I_\cC)$ 
of the configuration
space boundary an $(n-1)d-1$--dimensional unit sphere of velocities
$w = w^E_\cC +  w^I_\cC$. The internal part is split further into the component
$\langle w^I_\cC,Q^I_\cC\rangle Q^I_\cC$ parallel to
$Q^I_\cC$ and  the one perpendicular to it. In order to obtain a 
somewhat simpler form of the
differential equation, we rescale the parallel part, setting
\beq
v^I_\cC := \big(W_\cC(Q^I_\cC)\big)^{1/2} \, \langle w^I_\cC,Q^I_\cC\rangle \ \mbox{ and }\
X^I_\cC := w^I_\cC - \langle w^I_\cC, Q^I_\cC \rangle Q^I_\cC \, .
\Leq{def:vI:XI} 
\begin{remark}[gradient-like flow]\quad\\ 
We defined $v^I_\cC$ to be  $\langle w^I_\cC,Q^I_\cC\rangle$  rescaled  by 
$\sqrt{W_\cC(Q^I_\cC)}$ in order that $v^I_\cC$ becomes a  Lyapunov 
function on the collision manifold $\Xi_\cC$ when $w^E_\cC = 0$.
In this way the flow on the invariant manifold $\Xi_\cC \cap \{ w^E_\cC = 0 \}$
is  gradient-like\footnote{{\bf Definition \cite{McG1}: }
Let $\psi$ be a flow on a complete metric space $X$. Suppose
there is a continuous function $g: X \to \bR$ such that $g(\psi(x,t)) < g(x)$
if $t>0$
unless $x$ is a rest point. Suppose further that the rest points of $\psi$ are isolated. 
Then $g$ is called {\bf gradient-like}.\\
We use this definition, although in our case the rest points are not isolated.}. 
Condition  $w^E_\cC = 0$  is important since  orbits in $\widehat{\Sigma}_E$ colliding with $\Xi_\cC$ 
must satisfy this condition in the limit:  
\beq 
\lim_{\tau\to\infty} w^E_\cC(\tau) = 0 \, .
\Leq{limit:w_E}
In order to derive equation (\ref{limit:w_E}) we use that 
the limiting internal cluster energies~\footnote{Here we use a simplified notation 
for the values of observables along orbits.}
 $\lim_{\tau\to\infty} H^I_C(\tau)$  (see \eqref{HcC:I}) exist and are finite for all  
  $C\in \cC$. (Their existence  has been proven for a more general class of potentials
in \cite[Corollary 5.7]{FK}.) The existence of  $\lim_{\tau\to\infty} U^E_\cC(\tau)\in \bR$
for the external cluster energy is obvious. 
 $\lim_{\tau\to\infty} U^I_C(\tau) = + \infty$  so the finiteness of the limit of $H^I_C(\tau)$ implies that   
$\lim_{\tau\to\infty} K^I_C(\tau)  = + \infty$.
On the other hand  $\lim_{\tau\to\infty} K^E_C(\tau)\in \bR$, since the cluster-external 
forces are bounded along the orbit. Consequently, viewed projectively, the ratio of  
 internal and external speeds 
({\em i.e.} of $\big[ \sqrt{2K^I_C},  \sqrt{2K^E_C} \big]$)  tends to $[1, 0]$.
Now   \eqref{limit:w_E} follows from this fact,  \eqref{UUU} and
the definition \eqref{w:tau} of $w = w^E_\cC +  w^I_\cC$.
\hfill $\Diamond$
\end{remark}
\begin{lemma}[Dynamics at collisions] \quad \label{lem:diff:equations:local}\\
For the time variable $\tau$, see \eqref{new:time}, and the coordinates
\eqref{def:r:Q} and \eqref{def:vI:XI} 
the restriction of the o.d.e.\ \eqref{qw:p}
to the $\cC$ component of the boundary is of the form
\begin{align}
(r^I_\cC)' = &\, 0 
\qmbox{,}
(Q^I_\cC)' = \big(W_\cC(Q^I_\cC) \big)^{-1/\alpha} \, X^I_\cC \, ,
\label{rI:QI:p}\\
(v^I_\cC)' =&\, 
\big(W_\cC(Q^I_\cC)\big)^{1/2-1/\alpha}\, 
\Big[  (1- {\textstyle \frac{\alpha}{2}}) \,\big(\|w^I_\cC\|^2-\langle w^I_\cC,Q^I_\cC\rangle^2\big) 
- {\textstyle \frac{\alpha}{2}} \|w^E_\cC\|^2 \Big] \, ,
 \label{vI:p}\\[2mm]
(X^I_\cC)' =&\, 
F^I_\cC(Q^I_\cC) - \langle F^I_\cC(Q^I_\cC)  , Q^I_\cC \rangle Q^I_\cC
-   \langle F^I_\cC(Q^I_\cC) + W_\cC(Q^I_\cC)^{-1/\alpha}Q^I_\cC, w^I_\cC \rangle X^I_\cC
\NN \\
&- W_\cC(Q^I_\cC)^{-1/\alpha} \langle X^I_\cC, w^I_\cC \rangle Q^I_\cC 
\label{XI:p} \, ,\\[2mm]
(q^E_\cC)' =&\, 0 \qmbox{and} 
(w^E_\cC)'= - \langle F^I_\cC(Q^I_\cC),w^I_\cC \rangle \, w^E_\cC \, .
\label{qwE:p}
\end{align}
In particular,  the boundary component $r=0$,
is flow-invariant and so is  its submanifold $r = 0 = w^E_\cC$. Finally 
$v^I_\cC$ is strictly increasing on this submanifold at points at which   
$w^I_\cC$ and $Q^I_\cC$ are linear independent. 
\end{lemma}
\textbf{Proof:}
We derive these identities from the differential equation \eqref{qw:p}, 
beginning with \eqref{rI:QI:p}. For $q = q_0 + r Q^I_\cC\in \widehat{M}$
\beq
r' = \frac{\langle q^I_\cC, (q^I_\cC)'\rangle }{r} \stackrel{\eqref{qw:p}}{=}
 \widetilde{G}_E(q)\,  \langle Q^I_\cC, w^I_\cC\rangle \, .  
\Leq{r:p}
As $r\searrow 0$, $ \widetilde{G}_E( q_0 + r Q^I_\cC)\searrow 0$ (asymptotically linearly),
whereas 
$\|Q^I_\cC\|_\cM = 1$ and  $\|w^I_\cC\|_\cM \le 1$, proving $r'=0$ on the boundary. Similarly,
\[ (Q^I_\cC)' = \Big(\frac{q^I_\cC}{r}\Big)' 
\stackrel{\eqref{r:p}}{=} 
 \frac{\widetilde{G}_E(q)}{r}\, \big(w^I_\cC - \langle Q^I_\cC, w^I_\cC\rangle Q^I_\cC \big)
 \ \stackrel{r\searrow 0}{\longrightarrow} \
\frac{w^I_\cC - \langle Q^I_\cC, w^I_\cC\rangle Q^I_\cC}
{\big(W_\cC(Q^I_\cC)\big)^{\frac{1}{\alpha}}} \,. \]
To derive \eqref{vI:p}, based on Definition \eqref{def:vI:XI} we write the derivative of
$v^I_\cC$ as the sum $(v^I_\cC)' = I + I\!I + I\!I\!I $ with 
\[I :=  \big(\big(W_\cC(Q^I_\cC)\big)^{1/2}\big)' \, \langle w^I_\cC,Q^I_\cC\rangle \qmbox{and}
I\!I :=  \big(W_\cC(Q^I_\cC)\big)^{1/2} \, \langle (w^I_\cC)',Q^I_\cC\rangle \, .\]
Using \eqref{W:C:p} and \eqref{rI:QI:p},
\[ I =  \textstyle{\frac12}\big(W_\cC(Q^I_\cC)\big)^{-1/2-1/\alpha} 
\Big\langle  \nabla U^I_\cC(Q^I_\cC)  -
\langle \nabla U^I_\cC(Q^I_\cC), Q^I_\cC\rangle\, Q^I_\cC\,,\, w^I_\cC\Big\rangle \, 
\langle w^I_\cC,Q^I_\cC\rangle \, , \]
whereas by \eqref{qw:p} and the definition \eqref{force:term} of $F$
\[ I\!I =\textstyle{\frac12} \big(W_\cC(Q^I_\cC)\big)^{-1/2-1/\alpha} 
\Big\langle  \nabla U^I_\cC(Q^I_\cC) - \langle  \nabla U^I_\cC(Q^I_\cC), w^I_\cC\rangle  w^I_\cC\, ,\, Q^I_\cC\Big\rangle \,.\]
So by ($-\alpha$)--homogeneity of $U^I_\cC$
\begin{align*}
I + I\!I & = - {\textstyle \frac{\alpha}{2}}\, \big(W_\cC(Q^I_\cC)\big)^{+1/2-1/\alpha} 
\big( 1-\langle w^I_\cC\, ,\, Q^I_\cC\rangle^2 \big) \,.
\end{align*}
In this expression we substitute $1 = \|w\|^2 = \|w^I_\cC\|^2 + \|w^E_\cC\|^2$.
Finally, by \eqref{rI:QI:p} 
\[I\!I\!I =  \big(W_\cC(Q^I_\cC)\big)^{1/2} \, \langle w^I_\cC,(Q^I_\cC)' \rangle =
 \big(W_\cC(Q^I_\cC)\big)^{+1/2-1/\alpha} 
\big(\|w^I_\cC\|^2-\langle w^I_\cC\, ,\, Q^I_\cC\rangle^2 \big) \, , \]
proving \eqref{vI:p}.\\[2mm]
The proof of \eqref{XI:p} uses \eqref{rI:QI:p} and 
$(w^I_\cC)' = F^I_\cC(Q^I_\cC) - \langle F^I_\cC(Q^I_\cC) , w^I_\cC \rangle w^I_\cC$.\\[2mm]
In \eqref{qwE:p} the equation  $(q^E_\cC)' =0$ 
follows, since $(q^E_\cC)' =\widetilde{G}_E(q) w^E_\cC$, with the
scaling factor $\widetilde{G}_E(q)$ being zero 
at the boundary.\\
The equation  for $(w^E_\cC)' $ is a consequence of 
$(w^E_\cC)' = \Pi^E_\cC F(q) - \langle F(q),w \rangle \, w^E_\cC$,
where the velocity $w$  has norm one
and the force term $F(q)$ goes to
$F(Q^I_\cC)$  for $q=q_0+r Q^I_\cC$ in the limit $r\searrow 0$,
whereas $ \Pi^E_\cC F(q)\to 0$ for the external force.
\hfill $\Box$
\subsection{Blowing up the energy surface at infinity}\label{sub:energy:infty}
%
The configuration sphere $\bS$ at  spatial infinity is the disjoint union of 
$\widehat{\bS} = \bS\setminus \Delta$ and $\bS\cap \Delta$, 
where the last set corresponds
to non-trivial clusters whose barycenters tend to infinity. 
We treat  $\widehat{\bS}$ first.     
\subsubsection{The case of single particles}
%
We first consider the motion of  $n$ single particles approaching $\widehat{\bS}$, 
which in our center of mass configuration space can happen for energies $E\ge0$ 
if $n\ge 2$. For $n=2$ we have $\widehat{\bS} = \bS$.
In the Figure \ref{fig:blownupConfigSpace} for three particles on the line, 
$\widehat{\bS}$ corresponds to the six segments of the outer circle.\\
Instead of $q = r \, Q\in \widehat M$ we use the polar coordinates $(z,Q)$ near infinity, with
\beq
z:= r^{-\alpha}\qmbox{and}\overline{q} := z\, Q = r^{-\alpha} \, Q= r^{-\alpha-1}q \, .
\Leq{polar:coordinates} 
The configuration space is mirror-symmetric:
$\widehat{M} =\{ (z,Q) \mid q\in \widehat{M} \}$.

In these coordinates the sphere at infinity corresponds to $z =0$. We set
\beq 
\overline{G}_E(\overline{q}) :=  z^{1+1/\alpha} \, \widetilde{G}_E(q)
=  z^{1+1/\alpha} \, \widetilde{G}_E(z^{-1/\alpha}Q)
\textstyle =  z \,\frac{E+ z U(Q)}{1+ z U(Q)} U^{-1/\alpha}(Q)
\Leq{G:bar}
(compare with \eqref{new:time}) and rewrite the force term \eqref{force:term}:
\[ \overline{F}(\overline{q}) := F(q) = \textstyle
z \frac{{\cal M}^{-1}\nabla U(Q)} {2(1 + z U(Q))\,U(Q)^{\frac{1}{\alpha}}} \, .\] 
The differential equation \eqref{qw:p} then takes the form
\beq
\overline{q}' =  \overline{G}_E(\overline{q}) \, \big(w-(1+\alpha)\langle w , Q \rangle \, Q \big) 
\qmbox{,}w' =    \overline{F}(\overline{q}) - \langle  \overline{F}(\overline{q}),w\rangle\, w  
\Leq{ode:at:infinity}
and in polar coordinates $(z,Q)$ one gets
\[ z' = - \overline{G}_E(\overline{q}) \,\langle Q,w \rangle\qmbox{,}
Q' = \textstyle{ \frac{\overline{G}_E(\overline{q})}{z}} \,(w - \langle Q , w \rangle \, Q) \, .\]
By looking at $\overline{G}_E$ in \eqref{G:bar}, we conclude that the right hand sides
of these differential equations are real-analytic for $Q\in\widehat{\bS}$ and $z\in[0,\infty)$.
At spatial infinity, that is, at $z=0$, they reduce to 
\[ z' = 0\qmbox{,}
 Q' = E\,  U^{-1/\alpha}(Q) \, (w - \langle Q , w \rangle \, Q) \qmbox{,} w'=0 \, .\]
In particular the boundary component of $\wideparen{\Sigma}_E$ over $\widehat{\bS}$
is invariant under the flow.
Whereas $Q' = 0$ for $E=0$, for $E>0$ this is the case if and only if $w\in\{-Q,Q\}$.
These velocities $w$ then correspond to the negative/positive time asymptotics of 
non-clustering particles.\\
Similar to the case of two bodies treated in Subsection \ref{sub:two:bodies}, the 
unparametrized motion takes place on the invariant 
great circle $\bS\cap\mathrm{span}(Q,w)$ in configuration space. 
The reparameterization is given by integrating the factor $E\,  U^{-1/\alpha}(Q)$. 
\begin{enumerate}[$\bullet$]
\item 
If these do not meet $\Delta$,
then the span is a half-circle.
As in Subsection \ref{sub:two:bodies}, 
we cannot connect asymptotically free solutions of the 
original Hamiltonian differential equation \eqref{Ham:eq} via the
flow at infinity. 
The reason is that (unlike the collision orbits and the brake orbits) they 
converge to the {\em stable} manifold. 
That is nice, because this would not have a sensible physical interpretation
(observe that then the original time variable $t$ diverges as $\tau$ does).
\item 
Otherwise the trajectory can be asymptotic to points
$Q\in \Delta\cap \bS$, where $U^{-1/\alpha}(Q)$ vanishes.
This leads us to the next point, the blow-up for non-trivial clusters at spatial infinity.
\end{enumerate}
\subsubsection{The case of non-trivial clusters}
%
We now consider trajectories in $M$ that approach a point $q_0\in \Delta\cap \bS$.
In the Figure \ref{fig:blownupConfigSpace} for three particles on the line, 
this part of the boundary of blown up configuaration space $\wideparen{M}$ 
corresponds to the twelve small quarter circle segments near the the outer circle.\\
By the stratification \eqref{Delta:decomposition} of $\Delta$, 
$q_0 \in \Xi_\cC \cap \bS$ for a unique $\cC\in \cP_\Delta(N)$, 
with $\cC \neq \cC_{\max}$, since total collision occurs at $0\in M$.
$\Xi_\cC$, defined in \eqref{Xi:Null}, 
is a relatively open (and dense) subset of the linear subspace $\Delta^E_\cC$.
So by \eqref{dim:MC} it is a manifold of dimension $d (|\cC| - 1)$, with a number
$2\leq |\cC| \leq n-1$ of clusters. It follows that  $\Xi_\cC\cap \bS$ is open in a 
$\big(d (|\cC| - 1) - 1\big)$--dimensional sub-sphere of $\bS$.
Its blow-up is a fiber bundle
\beq 
B_\cC\ \longrightarrow \ \Xi_\cC\cap \bS \ 
\mbox{ , with typical fiber } S^{d (n - |\cC|) }_+ \, ,
\Leq{half-sphere:bundle} 
the plus sign meaning the half-sphere of incoming directions.\footnote{In general,
\eqref{half-sphere:bundle} is a  non-trivial bundle.
However, here we argue only semi-locally so that this does not play a role in our analysis.}
So it is a manifold whose
dimension coincides with the one of $\bS$, as it should, being part of 
$\partial \wideparen{M}$.

Next we construct coordinates on the half-sphere $S^{d (n - |\cC|) }_+$at $q_0$. 
The important result will be that, with the Cartesian coordinates 
$(q^E_\cC,q^I_\cC)$ near $\Xi_\cC$, $q^I_\cC$ can be used. 

We use the polar coordinates $(z,Q)$ of $\overline{M} = M\sqcup \bS$ near $\bS$ 
with $z := 1/r > 0$ for $q = r \, Q\in  M\backslash\{0\})$, and $z=0$ on $\bS$.
The metric used is
\[\big\|(z_1,Q_1) - (z_2,Q_2)\big\| := \sqrt{|z_1-z_2|^2 +  \|Q_1-Q_2\|^2} \ \quad
\big((z_j,Q_j)\in [0,\infty)\times \bS \big) .\]
For a unit vector $Q^I_\cC \in S^I_\cC$
we set $(q^E_\cC,q^I_\cC) := (c\, q_0, \, d\, Q^I_\cC)$ and consider the distance of
\[(z_1,Q_1) := \textstyle
\Big(\| q^E_\cC + q^I_\cC \|^{-1}\, , \, 
\frac{q^E_\cC + q^I_\cC}{\| q^E_\cC + q^I_\cC \|} \Big) \qmbox{and}
(z_2,Q_2) := (0,q_0) \, .\]
In the limit $c \nearrow \infty$, with $d\ge0$ fixed their difference is asymptotic to
\[(z_1,Q_1)-(z_2,Q_2)  = \textstyle
 \big( c^2+d^2)^{-1/2}, \frac{cq_0+dQ^I_\cC}{\sqrt{c^2+d^2}}-q_0 \big)\ 
\sim (\textstyle{\frac{1}{c},\frac{d}{c}} \,Q^I_\cC ) \, ,\]
which we write as $(\textstyle{\frac{1}{c},\frac{d}{c}}\,Q^I_\cC) 
= R\, \big( \cos(\vv), \sin(\vv)\,Q^I_\cC \big)$ 
(so that $\|(\textstyle{\frac{1}{c},\frac{d}{c}}\,Q^I_\cC)\|=R$ 
with $R=\sqrt{1+d^2}/c\searrow 0$).
As
\[q^I_\cC = d\, Q^I_\cC = \tan(\vv)Q^I_\cC \, ,\]
in the limit $d\nearrow \infty$ the $\cC$--internal cluster coordinates $q^I_\cC$ 
parameterize the half-sphere at $q_0$ of the bundle \eqref{half-sphere:bundle}. 
In that limit the distance between the positions of any pair of particles in different clusters
of $\cC$ goes to infinity, whereas the difference of particle positions in the same cluster
is constant.\\ 
From \eqref{qw:p} we get in this limit and for the cluster-internal potential
$U^I_\cC$ (see \eqref{VcC:I}),
the differential equation
\beq
(q^I_\cC)' =  \widetilde{G}_{E,\cC}(q^I_\cC) \, w^I_\cC \qmbox{,} 
(w^I_\cC)' =   F_\cC(q^I_\cC) - \langle F_\cC(q^I_\cC),w^I_\cC\rangle\, w^I_\cC,  
\Leq{qw:p:C}
with 
\[\textstyle
\widetilde{G}_{E,\cC}(q^I_\cC) := 
\frac{E+U^I_\cC(q^I_\cC)}
{( 1+U^I_\cC(q^I_\cC))\,(U^I_\cC(q^I_\cC))^{\frac{1}{\alpha}}} 
\qmbox{and} F_\cC(q^I_\cC) :=  
\frac{{\cal M}^{-1}\nabla U^I_\cC(q^I_\cC)} {2(1+U^I_\cC(q^I_\cC))\,
(U^I_\cC(q^I_\cC))^{\frac{1}{\alpha}}} \,.\]
So up to a common time reparameterization by 
$1/\big((1+U^I_\cC)((U^I_\cC)^{\frac{1}{\alpha}}\big)$, 
the clusters in $\cC$ only interact internally. To see this, note that, up to that  
factor, $\widetilde{G}_{E,\cC}$
and $F_\cC$ depend affinely on $U^I_\cC$.

The motion of the cluster centers occurs with a velocity vector $w^E_\cC$ 
that is a constant of the motion.

Like in Subsections \ref{sub:energy:coll} and the present one, 
with the half-sphere bundle $B_\cC$ from \eqref{half-sphere:bundle},
the boundary component
$B_\cC\times S^{(n-1)d-1}$ of the blown-up energy surface $\wideparen{\Sigma}_E$
is invariant under the flow, 
which by the same arguments is continuously differentiable.\\ 
$\cC$-internal collisions  that would lead to
components with a $\cD\in \cP(N)$ with $\cD$ strictly coarser than $\cC$ only can 
occur in the temporal limits $\tau\to\pm\infty$.\\
The same statement holds for $\|q^I_\cC(\tau)\|\to\infty$, 
that is escape to spatial infinity.\footnote{In Figure \ref{fig:blownupConfigSpace} for collisions these correspond  to the twelve 
dark green points near the outer circle. For  escape to spatial infinity, they are represented by
the twelve points on that circle.}
%
\section{Topology of the blown up configuration space}\label{sec:topology}
\begin{figure}[h]
\centerline{\includegraphics[width=68mm]{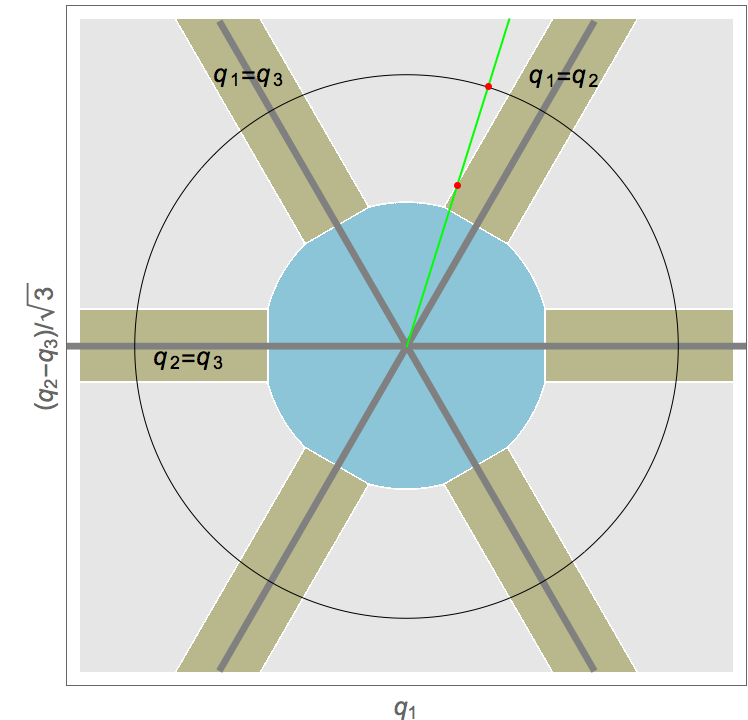}
\hfill
\includegraphics[width=66mm]{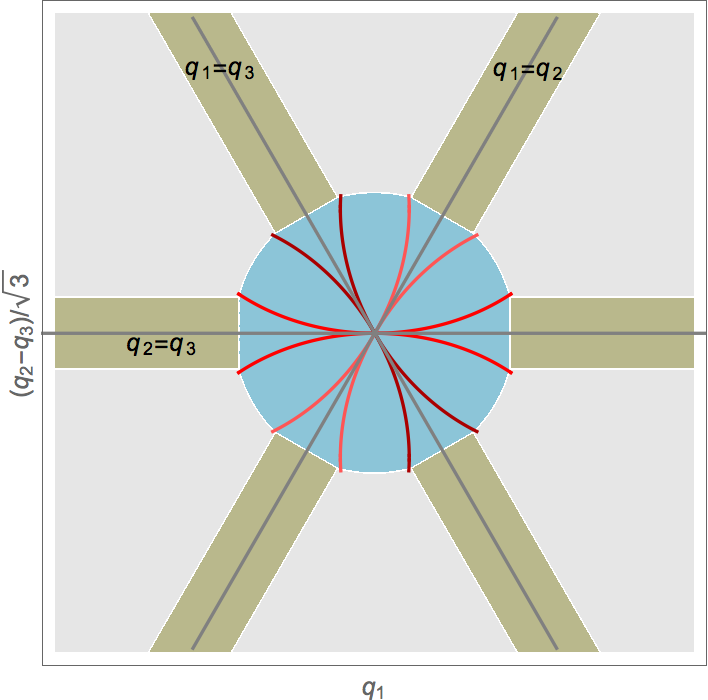}}
\caption{A Graf partition of the configuration space of $n=3$ particles in $d=1$ dimension, 
in center of mass system.
Left: the homeomorphism between $S^1\backslash \Delta$ and 
$\partial \Xi_{\cC_{\min}}^{(\delta_0)}$,
induced by the rays $R_c$ (green), with $\Xi_{\cC_{\min}}^{(\delta_0)}$ in gray.
Right: curves $S_\mathrm{ch}$, 
parametrised by $\delta<\delta_0$, corresponding to intersections of three atoms (red)} 
\label{fig:GrafBoundaryBlowUp}
\end{figure}

%
Here we are going to determine the homeomorphism type of the total boundary
blow up $\partial \wideparen{M}_\Delta$, with $\wideparen{M}_\Delta$ 
defined in Lemma \ref{lem:graph:blow:up}.
For the proof we use a variant of the Graf partition of configuration space $M$ devised by
{\sc Gian Michele Graf} in \cite{Gr}, see also \cite[Sect.\ 5.2]{DG}.
\begin{defi}{\bf (\cite[Sect.\ 12.6]{Kn})}\label{Graf--Partition}
For $\delta\in(0,1)$, let 
\[
J^{(\delta)}:M\to\bR\qmbox{,} J^{(\delta)}(q)
:= \max \big\{ J_\cC^E(q)+\delta^{|\cC|} \;\big|\;\cC\in\cP(N)
\big\}\, . 
\]
The \textbf{Graf partition} of the configuration space  $M$ is the
family of subsets 
\beq
\Xi_\cC^{(\delta)}:=
\l\{q\in M\;\Big|\; J_\cC^E(q) + \delta^{|\cC|} = J^{(\delta)}(q) \ri\}
\qquad \bigl(\cC\in\cP(N)\bigr).
\Leq{Z4}
\end{defi}
\begin{remark}[Graf partition] \quad\\[-6mm]\label{rem:graf}
\begin{enumerate}[1.]
\item 
These atoms are closed, and we obtain a measure theoretic partition of $M$: 
$\bigcup_{\cC\in\cP(N)}\Xi_\cC^{(\delta)}=M$,
For $\cC\neq\cD$ 
the Lebesgue measure of $\Xi_\cC^{(\delta)}\cap\Xi_\cD^{(\delta)}$
is zero, since the values of $J_\cC^E+\delta^{|\cC|}$ and 
$J_\cD^E+\delta^{|\cD|}$ coincide only on quadrics in~$M$.
\item 
Moreover, there is a $\delta_0\in (0,1)$ so that for all  $\delta\in(0,\delta_0]$, 
the  Graf partition (\ref{Z4})
has the property that for  $\Xi_\cC^{(\delta)}\cap\Xi_\cD^{(\delta)}\ne\es$,
the cluster decompositions  $\cC$ and $\cD$ are \textbf{comparable},
{\em i.e.}, $\cC\preccurlyeq\cD$ or $\cC\succcurlyeq\cD$ (see \cite[Lemma 12.52]{Kn}).
\hfill $\Diamond$
\end{enumerate}
\end{remark}
\begin{theorem}  \label{thm:boundary:hom}
For $n\ge 2$ particles, 
\begin{enumerate}[{\bf 1)}]
\item 
$S^{(n-1)d-1}\setminus \Delta$ is homeomorphic to $\partial\,\Xi_{\cC_{\min}}^{(\delta)}$,
\item 
whereas $\partial\,\Xi_{\cC_{\min}}^{(\delta)}$ is homeomorphic to 
$\partial \wideparen{M}_\Delta$.
\end{enumerate}
\end{theorem}
\textbf{Proof:}\\
{\bf 1)}
We first construct homeomorphisms
\[ H^{(\delta)}: S^{(n-1)d-1}\setminus \Delta\ \to\ \partial\,\Xi_{\cC_{\min}}^{(\delta)}
\qquad \big(\delta\in(0,\delta_0) \big)\]
to the boundary of the free atom $\Xi_{\cC_{\min}}^{(\delta)}$.
With the rays 
\[R_s:=\{\lambda s\mid \lambda>0\}  \ \subseteq \ M \qquad(s\in S^{nd-1}\setminus \Delta)\] 
$H^{(\delta)}(s)$ is defined as the unique intersection point in
$R_s\cap \partial\,\Xi_{\cC_{\min}}^{(\delta)}$.
In fact, 
\begin{enumerate}[1.]
\item 
$R_s\cap \partial\,\Xi_{\cC_{\min}}^{(\delta)}$ is non-empty:
\begin{enumerate}[(a)]
\item 
As $n\ge2$, there exists a set partition $\cD \in \cP_\Delta(N)$, and the difference 
\beq   \big[J_{\cC_{\min}}^E + \delta^n \big] - \big[J_\cD^E + \delta^{|\cD|} \big] 
= J_\cD^I + \delta^n - \delta^{|\cD|} 
\Leq{differ}
is strictly increasing along the ray $R_s$ and goes to $+\infty$, since the continuous map
$\lambda\mapsto J_\cD^I(\lambda s)$ equals $j\lambda^2$, with  $j=J_\cD^I(s) > 0$ if
$s\in S^{(n-1)d-1}\setminus \Delta$ (and $j = 0$ in the excluded case 
$s\in S^{(n-1)d-1}\cap \Delta_\cD^E\subseteq S^{(n-1)d-1}\cap \Delta$).
\item 
Conversely, 
$\lim_{\lambda\searrow 0} J_\cD^I(\lambda s) + \delta^n - \delta^{|\cD|} = \delta^n - \delta^{|\cD|}<0$ 
in \eqref{differ}, since $|\cD|<n$ and $\delta\in(0,1)$. 
\end{enumerate}
\item 
$R_s\cap \partial\,\Xi_{\cC_{\min}}^{(\delta)}$ consists of one point only, 
and the intersection is transverse. 
This follows since for $\lambda_0 s\in \Xi_{\cC_{\min}}^{(\delta)} \cap \Xi_{\cD}^{(\delta)}$,
\eqref{differ} is strictly increasing at $\lambda_0 s$, so that 
$\lambda s\in \Xi_{\cC_{\min}}^{(\delta)}$
for $\lambda > \lambda_0$. 
Of course, $\lambda_0 s$ can belong to more than two atoms in the
partition, see Figure \ref{fig:GrafBoundaryBlowUp}. Then that boundary point lies in the intersection
of more than one quadric, but $R_s$ is transversality to all of them.  
\item 
By that transversality property, $H^{(\delta)}$ is continuous. 
$H^{(\delta)}$  is injective, since the rays $R_s$ and 
$R_{s'}$ are mutually disjoint for $s\neq s'\in S^{(n-1)d-1}$. \\ 
As the atom $\Xi_{\cC_{\min}}^{(\delta)}$ is closed and disjoint from $\Delta$, the same is 
true for its boundary $\partial\,\Xi_{\cC_{\min}}^{(\delta)}$.
Since 
$\bigcup_{s\in S^{(n-1)d-1}\setminus \Delta} R_s
= \widehat{M}=\bR^{(n-1)d}\setminus \Delta$, 
$H^{(\delta)}$ is  surjective. \\
The inverse of $H^{(\delta)}$ is  continuous, too, since the intersection of the rays 
$R_s$ with the sphere is transverse, too.
\end{enumerate}
Incidentally, we proved that $H^{(\delta)}$ is even locally bi-Lipschitz.\\[2mm]
So for $\delta_1<\delta_2\in(0,\delta_0)$ the boundaries 
$\partial\,\Xi_{\cC_{\min}}^{(\delta_1)}$
and $\partial\,\Xi_{\cC_{\min}}^{(\delta_2)}$ are homeomorphic. Furthermore, 
${\rm int} \big(\, \Xi_{\cC_{\min}}^{(\delta_1)} \big) \supseteq \Xi_{\cC_{\min}}^{(\delta_2)}$ 
if $\delta_0\le1/2$, and 
$\lim_{\delta\searrow0} \Xi_{\cC_{\min}}^{(\delta)}= \widehat{M}$ 
in the sense of Hausdorff distance, 
since $0 < |\cD|<n$ for $\cD \in \cP_\Delta(N)$ so that in \eqref{differ}
\beq
 0 < \delta_1^{|\cD|} - \delta_1^n < \delta_2^{|\cD|} - \delta_2^n\qmbox{and}
\lim_{\delta\searrow0} \big(\delta^{|\cD|} - \delta^n \big) = 0.
\Leq{ineq:dellta}
\bigskip

\noindent
{\bf 2)}
To relate the boundaries $\partial\,\Xi_{\cC_{\min}}^{(\delta)}$ and $\partial \wideparen{M}_\Delta$, we use
homeomorphisms different from 
$H^{(\delta_1)}\circ (H^{(\delta_2)})^{-1}: 
\partial\,\Xi_{\cC_{\min}}^{(\delta_2)}\to \partial\,\Xi_{\cC_{\min}}^{(\delta_1)} $, since the limit
of the latter for 
$\delta_1\searrow 0$ is not well-behaved. Instead we show existence of a 
smooth vector field
\beq
v: \widehat{M} \setminus \Xi_{\cC_{\min}}^{(\delta_0)}\to \bR^{d(n-1)},
\Leq{vector:field} 
whose time $t$ flow restricts to locally Lipschitz homeomorphisms 
\beq
\partial\,\Xi_{\cC_{\min}}^{(\delta_2)}\to \partial\,\Xi_{\cC_{\min}}^{(\delta_1)} 
\qquad \big(\delta_1\in (0,\delta_2]\mbox{ with }\delta_2-\delta_1=t\big).
\Leq{time:t:flow} 
Every point $q\in\widehat{M} \setminus \Xi_{\cC_{\min}}^{(\delta_0)}$ belongs to
 $q\in\partial\,\Xi_{\cC_{\min}}^{(\delta)}$, for some $\delta\in(0,\delta_0)$.
  
Furthermore $q$ (exactly) belongs to
$\bigcap_{\ell = 0}^k \partial\,\Xi_{\cD_\ell}^{(\delta)}$ with $\cD_0 := \cC_{\min}$ and
$k\in\bN$.
According to Remark \ref{rem:graf}.2 these set partitions are mutually comparable. 
So by reordering we can assume that the chain
\beq
\mathrm{ch} := \{\cD_0,\ldots, \cD_k\}\mbox{ is ordered: }\cD_0\preccurlyeq \cD_1\preccurlyeq \ldots\preccurlyeq \cD_k,\quad(1\le k\le n-1).
\Leq{chain}
Conversely, every such chain \eqref{chain} defines a semialgebraic set $S_\mathrm{ch}\subseteq 
\widehat{M} \setminus \Xi_{\cC_{\min}}^{(\delta_0)}$, consisting of the $q$ that for some 
$\delta\in(0,\delta_0)$ simultaneously belong to all $\partial\,\Xi_{\cD_\ell}^{(\delta)}$, 
see Figure \ref{fig:GrafBoundaryBlowUp}, right. Finally, the set of all chains \eqref{chain} 
gives rise to a set partition of $\widehat{M} \setminus \Xi_{\cC_{\min}}^{(\delta_0)}$ 
by the $S_\mathrm{ch}$.
In fact, these semialgebraic sets are submanifolds, of codimension $k-1$.

For a chain \eqref{chain} the
corresponding projections \eqref{proj:E} are pairwise commuting.
So the level sets $(J_{\cD_\ell}^I - J_{\cD_{\ell-1}}^I)^{-1}(\delta^{|\cD_\ell|} - \delta^{|\cD_{\ell-1}|})$ 
are pairwise $\cM$-orthogonal
and define smooth functions 
\[ \tilde{f}_\ell: U_\mathrm{ch} \longrightarrow \bR^+ \ \mbox{ ,}\quad 
\big( J_{\cD_\ell}^I - J_{\cD_{\ell-1}}^I \big)^{-1} \big(\delta^{|\cD_\ell|} - \delta^{|\cD_{\ell-1}|} \big) 
\ni q\longmapsto \delta 
\qquad(\ell=1,\ldots, k)\]
on a suitable neighborhood
$U_\mathrm{ch}\subseteq \widehat{M} \setminus \Xi_{\cC_{\min}}^{(\delta_0)}$ of $S_\mathrm{ch}$ 
that have non-vanishing gradients. 
Formulated differently, the vector fields
\[ \tilde{v}_\ell: U_\mathrm{ch}\to\bR^{(n-1)d}\qmbox{,}
\tilde{v}_\ell := \textstyle 
\frac{- \nabla \tilde{f}_\ell} {\| \nabla \tilde{f}_\ell \|_\cM^2} \qquad (\ell = 1, \ldots , k)\]
are pairwise $\cM$-orthogonal, with Lie derivatives $L_{\tilde{v}_i} \tilde{f}_\ell = - \delta_{i,\ell}$.
Therefore there exists a unique linear combination 
\[v_\mathrm{ch} := \sum_{\ell=1}^k c_\ell\, \tilde{v}_\ell :U_\mathrm{ch}\to\bR^{(n-1)d}\]
of the vector fields so that $L_{v_\mathrm{ch} } f_\ell = - \delta_{i,\ell}$ for the functions
\beq
f_\ell : U_\mathrm{ch} \longrightarrow \bR^+ \qmbox{,} 
\big( J_{\cD_\ell}^I \big)^{-1} \big( \delta^{|\cD_\ell|} - \delta^n  \big) 
\ni q\longmapsto \delta  \qquad(\ell=1,\ldots, k) 
\Leq{f:ell}
that satisfy
\[f_\ell \big( \partial\,\Xi_{\cD_\ell}^{(\delta)} \cap U_\mathrm{ch} \big) = \delta\]
(note that $J_{\cD_0}^I = J_{\cC_{\min}}^I = 0$, so that 
$J_{\cD_\ell}^I = \sum_{m=1}^\ell (J_{\cD_m}^I - J_{\cD_{m-1}}^I) $).
In particular, $v_\mathrm{ch}$ is tangential to $S_\mathrm{ch}$.
Actually we define the neighborhood $U_\mathrm{ch} $ by the condition that for some constant $c>1$
\beq
\textstyle \frac{f_i(q)}{f_j(q)}\in (1/c,c)\qquad (i,j=1,\ldots,k).
\Leq{fi:fj}
Note this still guarantees that $S_\mathrm{ch}\subseteq U_\mathrm{ch}$.
Of course for $k=1$, the restriction \eqref{fi:fj} is vacuous.

$c>1$ is chosen small enough so that $U_\mathrm{ch} \cap U_{\mathrm{ch}'}=\emptyset$ if
the chains $\mathrm{ch}$ and $\mathrm{ch}'$ contain incompatible atoms.

Combining the vector fields $v_\mathrm{ch}$ by a partition of unity subordinate to the 
$U_\mathrm{ch}$, we get a vector field \eqref{vector:field} on 
$\widehat{M}\setminus \Xi_{\cC_{\min}}^{(\delta_0)}$.
It induces a local flow $\Phi$,  mapping level surfaces to level surfaces.  

Restricting the local 
flow to initial conditions in $\partial\,\Xi_{\cC_{\min}}^{(\delta_2)}$, we
obtain the family of homeomorphisms \eqref{time:t:flow}.
Their limit for $\delta_1\searrow 0$ exists (see below) and 
leads to a homeomorphism with the boundary blow up $\partial \wideparen{M}_\Delta$.

Existence of the limit is seen as follows: 
\begin{enumerate}[(a)]
\item
Along a $\Phi$ trajectory at time $t\ge0$, \eqref{fi:fj} transforms into $\frac{f_i(q)-t}{f_j(q)-t}$.
So if $q\in U_\mathrm{ch}\setminus S_\mathrm{ch}$, then its trajectory will ultimately leave $U_\mathrm{ch}$.
This implies that near the escape time limit we can assume that $\Phi_t(q)\in S_{\mathrm{ch}'}$ for 
some subchain $\mathrm{ch}'$.
\item 
For $q\in S_\mathrm{ch}$ with $f_\ell(q)=\delta$ by definition  
$\lim_{t\nearrow\delta}f_\ell\circ \Phi_t(q)=0$ ($\ell=1,\ldots,k$). This implies that
$\lim_{t\nearrow\delta} \Phi_t(q)\in \Delta$ exists:
\begin{enumerate}[$\bullet$]
\item 
The external coordinates $\big(\Phi_t(q)\big)_{\cD_k}^E$ are independent of $t$,
since \eqref{f:ell} depends only on internal $\cD_k$ coordinates, and by $\cM$-orthogonality
of the projections \eqref{proj:E}.
\item 
The internal coordinates $\big(\Phi_t(q)\big)_{\cD_k}^I$ go to zero as $t\nearrow\delta$.
\end{enumerate}
\item 
For the same reason, the map 
\beq
\Phi^{(\delta)}:
\partial\,\Xi_{\cC_{\min}}^{(\delta)}\to \Delta \qmbox{,} q\mapsto \lim_{t\nearrow\delta}\Phi_t(q)
\Leq{Phi:delta} 
is continuous.
It is surjective, since the Hausdorff distance of these two subsets of $M$ goes to zero as 
$\delta\searrow0$. It is not injective (not even in the case $n=2$, since then 
$\Delta=\Delta_{\cC_{\max}}^E =\{0\}$, whereas  
$\partial\,\Xi_{\cC_{\min}}^{(\delta)}\cong  S^{n-1}$).
\item 
We recall the definition $S^I_\cC=\big(J^I_\cC\big)^{-1}(1)$ of the internal unit sphere
for $\cC\in\cP(N)$ (which is of dimension $d(n-|\cC|)-1$).\\
For $q\in\partial\,\Xi_{\cC_{\min}}^{(\delta)} \cap \overline{S}_\mathrm{ch}$  
and $\ell=1,\ldots,k$, the unit vectors $N_{\cD_\ell}^{(\delta)}(q)$ given by
\beq
N_{\cD_\ell}^{(\delta)} : \partial\,\Xi_{\cC_{\min}}^{(\delta)} \cap \overline{S}_\mathrm{ch}
\to S^I_{\cD_\ell} \qmbox{,}
N_{\cD_\ell}^{(\delta)}(q) = \lim_{t\nearrow\delta} 
\frac{(\Phi_t(q) - \Phi^{(\delta)}(q))^I_{\cD_\ell}}
{\|(\Phi_t(q) - \Phi^{(\delta)}(q))^I_{\cD_\ell}\|}
\Leq{normal} 
exist, and obviously depend only on $\cD_\ell\in \cP_\Delta(N)$, not on the chain $\mathrm{ch}$
to which $\cD_\ell$ belongs.  This is important since the closure $\overline{S}_\mathrm{ch}$
of $S_\mathrm{ch}$ used in the definition of the domain can intersect other 
$S_{\mathrm{ch}'}$. Moreover, the maps $N_{\cD_\ell}^{(\delta)}$ are continuous, and 
$N_{\cD_\ell}^{(\delta)}(q)$ is perpendicular to $\Delta_{\cC_\ell}^E$.

By the above arguments, on the subsets 
$\partial\,\Xi_{\cC_{\min}}^{(\delta)} \cap S_\mathrm{ch}$, indexed by the $\mathrm{ch}$,
\item 
Together, $\Phi^{(\delta)}(q)$ and 
the $N_{\cD_\ell}^{(\delta)}(q)$ ($\ell=1,\ldots,k$) define a point in $\partial \wideparen{M}_\Delta$.
So \eqref{Phi:delta} and \eqref{normal} define a continuous map 
\[ \Psi^{(\delta)}: \partial\,\Xi_{\cC_{\min}}^{(\delta)}\to \partial \wideparen{M}_\Delta\]
The (continuous and surjective)
blow-down $\beta:\wideparen{M}_\Delta\to M$ maps $\wideparen{M}_\Delta$ onto $\Delta$, whereas it is the identity
on $\widehat{M}=M\setminus \Delta$. It has the property
\beq
\beta\circ \Psi^{(\delta)} = \Phi^{(\delta)} \, .
\Leq{blow:down}
\item 
We finally prove that $\Psi^{(\delta)}$ is a homeomorphism. 
\begin{enumerate}[$\bullet$]
\item 
To show that it is a surjection, we use that the right hand side of \eqref{blow:down} is a surjection.
For the chains $\mathrm{ch}=\{\cD_0,\cD_1\}$ of length $k=1$ the restriction of
$\Psi^{(\delta)}$ to $\partial\,\Xi_{\cC_{\min}}^{(\delta)}\cap S_\mathrm{ch}$ maps onto 
${\rm int}(\beta^{-1}(\Xi_{\cD_1}))\subseteq \partial \wideparen{M}_\Delta$, since then the vector field on 
$U_\mathrm{ch}$ is ultimately radial, so that the $t$-dependent unit vector in \eqref{normal} becomes 
constant. But the union of these interiors is dense in $\wideparen{M}_\Delta$, which, together with the
continuity of $\Psi^{(\delta)}$, shows that this map is onto. 
\item
To prove injectivity of $\Psi^{(\delta)}$, by the $\Phi$-invariance of the sets $S_\mathrm{ch}$ it
suffices to consider points $q_1,q_2 \in \Psi^{(\delta)} \cap S_\mathrm{ch}$ and conclude that
they coincide if their images coincide. 
Using (a) above, we can further assume that $v(q_i) = v_\mathrm{ch}(q_i)$. This property is then preserved by the forward flow.
By the first bullet point in (b), we can also assume (by diminishing $\delta>0$, if necessary) 
that their external $\cD_k$ coordinates coincide. Now if $q_1\neq q_2$, there is a largest 
$\ell\in \{1,\ldots,k\}$ such that $(q_1)_{\cD_\ell}^I\neq (q_2)_{\cD_\ell}^I$. The 
nontrivial rotation of the plane spanned by 
$(q_1)_{\cD_\ell}^I$ and $(q_2)_{\cD_\ell}^I\in \Delta_{\cD_\ell}^I$
mapping the first to the second point maps the vector field $v$ along the forward orbit 
$t\mapsto \Phi_t (q_1)$ onto the one of the forward orbit of $q_2$. Thus 
\[N_{\cD_\ell}^{(\delta)}(q_1)\neq N_{\cD_\ell}^{(\delta)}(q_2),\]
since (by diminishing $\delta>0$ again, if necessary), we can assume that 
$N_{\cD_\ell}^{(\delta)}(q_i)$ is not perpendicular to $q_{\cD_\ell}^{(\delta)}$.
\item 
Although the continuous bijection $\Psi^{(\delta)}$ does not have a compact domain, 
its inverse is continuous, too. Namely 
the intersections of  $\partial\,\Xi_{\cC_{\min}}^{(\delta)}$ with closed balls in $M$ of radius $r>0$ 
are compact, and exhaust the domain as $r\to\infty$
As the restrictions of $\Psi^{(\delta)}$ are homeomorphisms onto their images 
(since $\wideparen{M}_\Delta$ is Hausdorff),
this shows that $\Psi^{(\delta)}$ itself is a homeomorphism.~$\Box$
\end{enumerate}
\end{enumerate}

\begin{corollary}\quad\\
For $d=1$ dimensions, $\partial \wideparen{M}_\Delta$ 
is homeomorphic to $n!$ disjoint copies of $\bR	^{n-2}$. 
\end{corollary}
\textbf{Proof:}
This follows from Theorem \ref{thm:boundary:hom}, since for $d=1$ the set $\Delta$ 
is a union of
${n \choose 2}$ hyperplanes $\{q\in \bR^{n-1}\mid q_i=q_j\}$ ($1\le i<j\le n$), so that
$S^{n-2}\setminus \Delta$ is the disjoint union of $n!$ 
relatively open spherical simplices, given by
the ordering of the coordinates $(q_1,\ldots,q_n)$. 
These are in turn homeomorphic to $\bR^{n-2}$.
\hfill $\Box$
\appendix
\section{Appendix: Manifolds with corners}\label{sect:appendix}
%
We follow \cite{AMN} and \cite{Me} in our presentation. 
Manifolds with corners are modeled on the $m$--dimensional cylinders 
\[\bR^m_k := [0,\infty)^k\times \bR^{m-k}\,\subseteq \bR^m\qquad (k\le m\in\bN_0).\]
Their subsets
\beq
L_I := \{x = (x_1,\ldots,x_m)\in \bR^m_k\mid x_i=0\mbox{ if }i\in I \}\qquad(I\subseteq \{1,\ldots,m\})
\Leq{L:I}
will be used to define submanifolds.
\begin{defi}\label{def:smooth}
Let $U\subseteq \bR^m_k$ and $V\subseteq \bR^{m'}_{k'}$  be open, and $f:U\to V$.\quad\\[-6mm]
\begin{enumerate}[$\bullet$]
\item 
$f$ is called {\bf smooth} if for some open neighbourhood 
$\tilde{U}\subseteq \bR^m$ of $\,U$ there exists 
$\tilde{f}\in C^\infty\big(\tilde{U},\bR^{m'}\big)$ with $\tilde{f}|_U = f$.
\item 
$f$ is called a {\bf diffeomorphism}, if it is a smooth bijection with $f^{-1}$ smooth. 
\end{enumerate}
\end{defi}
\begin{defi}[manifolds with corners]
Let $X$ be a Hausdorff space.
\begin{enumerate}[$\bullet$]
\item 
An  {\bf ($m$--dimensional) corner chart} $(U,\phi)$ on $X$ is a homeomorphism 
$\phi:U\to V$, with $V$ open in $\bR^m_k$.
\item 
Corner charts $(U_1,\phi_1)$ and $(U_2,\phi_2)$ on $X$ are {\bf compatible} if for $U:=U_1\cap U_2$
\[\phi_2\circ \phi_1^{-1}:\phi_1(U)\to \phi_2(U)\]
is a diffeomorphism.
\item 
A {\bf (corner) atlas} $\{(U_i,\phi_i)\mid i\in I\}$ on $X$ is a family of pairwise compatible charts
$(U_i,\phi_i)$ on $X$ of equal dimension with $\bigcup_{i\in I}U_i=X$.
\item 
Corner atlases on $X$ are {\bf equivalent} if their union is a corner atlas on $X$.\linebreak
A {\bf corner structure} on $X$ is an equivalence class of corner atlases of $X$.
\item 
A paracompact Hausdorff space $X$ with a corner structure 
consisting of $m$--dimensional corner charts 
is an {\bf ($m$--dimensional) manifold with corners}.
\item 
For $\partial_\ell\bR^m_k:=
\{x\in \bR^m_k\mid \mbox{of } x_1,\ldots, x_k, \mbox{ exactly } \ell \mbox{ vanish}  \}$,
\[\partial_\ell X := \{p\in X\mid\mbox{coordinates at } p \mbox{ map to } \partial_\ell\bR^m_k\}\]
and the {\bf boundary} $\partial X:=  \partial^1X$ {\bf of} $X$ for 
$\partial^\ell X := \overline{\partial_\ell X}$.
\end{enumerate}
\end{defi}
Unlike for manifolds with boundary, the Cartesian product of two manifolds with corners is 
naturally a manifold with corners.
\begin{defi}[submanifolds of manifolds with corners]\label{def:sub}\quad\\[-5mm]
\begin{enumerate}[$\bullet$]
\item 
A subset $S\subseteq X$ of an $m$-dimensional manifold with corners is a 
{\bf weak submanifold} if for every $x\in S$ there exist $k\in\{1,\ldots,m\}$ and a 
corner chart $\phi: U\to\Omega\subseteq \bR^m_k$ with $x\in U$ such that $\phi(S\cap U)$ is a 
submanifold of $\,\bR^m$.\\
Then the {\bf dimension of $S$ at $x$} is $\dim(\phi(S\cap U))$ at $\phi(x)$.
\item 
A weak submanifold $S\subseteq X$ is a {\bf submanifold} (in the sense of manifolds with corners) if, 
additionally there are integers $m' \le m$ and $k'\le m'$, and a matrix
$G\in {\rm GL}(m, \bR)$ such that 
\begin{enumerate}
\item 
$G\cdot \big(\bR^{m'}_{k'}\times\{0\}\big) \subseteq \bR^m_k$
\item
The chart $\phi$ maps $S\cap U$ bijectively to the intersection of this linear submanifold
with $\Omega$, in other words
$\phi(S\cap U)=G\cdot \big(\bR^{m'}_{k'}\times\{0\}\big)\cap \Omega$.
\end{enumerate}
\item
A submanifold $S\subseteq X$ is a {\bf $p$--submanifold} if for $x\in X$ there exists a corner chart
$(U,\phi)$ at $x$ and $I\subseteq \{1,\ldots,m\}$ with, see Definition \eqref{L:I} 
\[\phi(S\cap U) = L_I\cap \phi(U).\]  
Then $|I|$ is the {\bf codimension of $S$ at $x$} and $|I\cap\{1,\ldots, k\}|$ is the 
{\bf boundary depths of $S$ at $x$}.
\end{enumerate}
\end{defi}
So a $p$-submanifold $S$ of $X$ is a closed submanifold that has a tubular neighborhood: 
$S\subseteq U\subseteq X$ that is locally of product form.
\addcontentsline{toc}{section}{References}

\begin{thebibliography}{99}

\bibitem[Ai]{Ai}
Martin Aigner: Combinatorial Theory. Classics in
Mathematics, Springer 1997

\bibitem[AMN]{AMN}
Bernd Ammann, J\'er\'emy Mougel, Victor Nistor:
A comparison of the Georgescu and Vasy spaces associated to the $n$-body problems.
Annales Henri Poincar\'e  {\bf 23}, 1141--1203 (2022)

\bibitem[De]{Dev}
Robert L. Devaney:
Singularities in classical mechanical systems. In: 
Ergodic theory and dynamical systems I. Birkh\"auser, Boston, MA, 1981

\bibitem[DG]{DG}
Jan Derezi\'{n}ski, Christian G\'{e}rard: 
Scattering Theory of Classical and Quantum $N$-Particle Systems. 
Texts and Monographs in Physics. Berlin: Springer 1997

\bibitem[DMMY]{DMMY}
Nathan Duignan, Richard Moeckel, Richard Montgomery, Guowei  Yu: 
Chazy-type asymptotics and hyperbolic scattering for the $n$-body problem. 
Archive for Rational Mechanics and Analysis {\bf 238}, 255--297  (2020)

\bibitem[ElB]{ElB}
Mohamed Sami ElBialy: Collision singularities in celestial mechanics.
SIAM journal on mathematical analysis {\bf 21}, 1563--1593 (1990)

\bibitem[FK]{FK}
Stefan Fleischer, Andreas Knauf: 
Improbability of Collisions in $n$-Body Systems. 
Archive for Rational Mechanics and Analysis {\bf 234}, 1007-1039, (2019)

\bibitem[FKM]{FKM}
Jacques F\'ejoz, Andreas Knauf, Richard Montgomery:
Lagrangian Relations and Linear Point Billiards,
Nonlinearity {\bf 30}, 1326--1355 (2017)

\bibitem[Gr]{Gr}
Gian Michele Graf: 
Asymptotic completeness for $N$-body short range quantum systems: A new proof.
Commun. Math. Phys. {\bf 132}, 73--101 (1990) 

\bibitem[He]{He}
Douglas Heggie: 
A global regularisation of the gravitational N-body problem. 
Celestial mechanics {\bf 10}, 217--241 (1974) 

\bibitem[Hi]{Hi}
Morris W. Hirsch: Differential Topology. 
Graduate Texts in Mathematics, Vol.\ 33. 
Berlin, Heidelberg, New York: Springer 1988

\bibitem[KK]{KK}
Andreas Knauf, Markus Krapf: 
The non-trapping degree of scattering. 
Nonlinearity {\bf 21}, 2023--2041 (2008)

\bibitem[KM]{KM}
Andreas Knauf, Richard Montgomery: 
Compactifying the three-body problem.
[In preparation]

\bibitem[Kn]{Kn}
Andreas Knauf: Mathematical Physics: Classical Mechanics.
Springer, New York, 2018.

\bibitem[LI]{LI}
Ernesto A. Lacomba, Luis A. Ibort:
Origin and infinity manifolds for mechanical systems with homogeneous potentials.
Acta Applicandae Mathematica {\bf 11}, 259--284 (1988)

\bibitem[LS]{LS}
Ernesto A. Lacomba, Carles Sim\'o: 
Boundary manifolds for energy surfaces in celestial mechanics.
Celestial Mechanics {\bf 28}, 37--48 (1982)

\bibitem[Le]{Le}  
Georges  Lema\^{i}tre: 
Regularization of the Three Body Problem.
Vistas in Astronomy  {\bf 1},  207--215  (1955)  

\bibitem[LC]{LC}  
Tullio Levi-Civita:  
Sur la r\'esolution qualitative du probl\`eme restreint. 
Acta Mathematica {\bf 30}, 305--327 (1906)

 \bibitem[MM]{MM}
John N. Mather, Richard McGehee:
Solutions of the collinear four body problem which become unbounded in finite time.
Dyn. Syst., Theor. Appl., Battelle Seattle 1974 Renc., 
Lect. Notes Phys. {\bf 38}, 573--597 (1975)

\bibitem[McG1]{McG1}
Richard McGehee:
Triple collision in the collinear three-body problem.
Invent. Math. {\bf 27}, 191--227 (1974)

\bibitem[McG2]{McG2}
Richard McGehee:
Double collisions for a classical particle system with nongravitational interactions.
Commentarii Mathematici Helvetici {\bf 56}, 524--557 (1981)

\bibitem[McG3]{McG3} 
Richard McGehee: 
A stable manifold theorem for degenerate fixed points with applications to celestial mechanics.
J. Differential Equations {\bf 14}, 70--88 (1973) 

\bibitem[Me]{Me}
Richard  Melrose: Differential analysis on manifolds with corners. 
Book in preparation (1996)

\bibitem[MM]{MM} 
Richard Moeckel, Richard Montgomery: 
Symmetric Regularization, Reduction, and Blow-Up of the Planar Three-Body Problem. 
Pac. J. Math. {\bf 262}, 129--189 (2013)

\bibitem[Mo]{Mo}
J\"urgen Moser: 
Regularization of Kepler's problem and the averaging method on a manifold.
CPAM {\bf 23}, 609--636  (1970)

\bibitem[Ro]{Ro} 
Clark Robinson: 
Homoclinic Orbits and Oscillations for the Planar Three-body Problem. 
J.\ Differential Equations {\bf 52},  356--377  (1984)

\bibitem[Va]{Va} 
Andr\'as Vasy: 
Propagation of singularities in three-body scattering.
Ast\'erisque {\bf 262}, 157 p. (2000) 

\bibitem[Wald]{Wald} 
J\"{o}rg Waldvogel: 
A new regularization of the problem of three bodies.
Celestial Mechanics {\bf 6},  221--231 (1972)

\bibitem[Wa]{Wa}
Qiu-Dong Wang: The global solution of the $n$--body problem.
Celestial Mechanics and Dynamical Astronomy {\bf 50}, 73--88 (1991)

\bibitem[Wall]{Wall} 
C.T.C. Wall:  Differential Topology, 
Cambridge studies in advanced mathematics, Cambridge University Press, 2016

\end{thebibliography}
\end{document}